\documentclass{article}

\usepackage{amsfonts,amsmath}
\usepackage{amssymb,latexsym}
\usepackage[dvips]{epsfig}

\title{\bf A Categorification of the Temperley-Lieb Algebra and 
Schur Quotients of $U(\mathfrak{sl}_2)$ via Projective and Zuckerman 
Functors\footnotetext{Published in 
Selecta Mathematica, New ser. 5 (1999), 199--241}}
\author{Joseph Bernstein, Igor Frenkel and Mikhail Khovanov}
\date{}
\newtheorem{prop}{Proposition}
\newtheorem{theorem}{Theorem}
\newtheorem{lemma}{Lemma}
\newtheorem{corollary}{Corollary}
\newtheorem{definition}{Definition}
\newtheorem{conjecture}{Conjecture}
\newcommand{\oplusop}[1]{{\mathop{\oplus}\limits_{#1}}}

\begin{document}
\maketitle
\baselineskip 12pt

\tableofcontents

\bigbreak

\def\C{\mathbb C}
\def\R{\mathbb R}
\def\N{\mathbb N}
\def\Z{\mathbb Z}
\def\Q{\mathbb Q}
\def\F{\mathbb F}
\def\Zq{\Z[q,q^{-1}]}
\def\l{\lbrace}
\def\r{\rbrace}
\def\o{\otimes}
\def\D{\Delta}
\def\O{\mathcal{O}}
\def\Ok{\mathcal{O}(k,n-k)}
\def\slt{\mathfrak{sl}_2}
\def\U{U_q(\mathfrak{sl}_2)}
\def\Ug{U_q(\mathfrak g)}
\def\fsl{\mathfrak{sl}_n}
\def\gl{\mathfrak{gl}_n}
\def\td{\textrm{dim}}
\def\g{\mathfrak{g}}
\def\sln{\mathfrak{sl}_n}
\def\p{\mathfrak{p}}
\def\CU{\mathcal{U}}
\def\opsi{\overline{\psi}}
\def\ophi{\overline{\phi}}

\def\lra{\longrightarrow}

\def\E{\mathcal E} 
\def\F{\mathcal F}

\def\dU{\dot{U}(\mathfrak{sl}_2)}
\def\dUq{\dot{U}_q(\mathfrak{sl}_2)}
\def\dB{\dot{\mathbb B}}

\def\UL{\mathbf U}
\def\UZ{\mathbf U_{\Z}}

\def\M{\mathcal M} 
\def\Mi#1{\M_{Zf}(#1 )}

\def\binom#1#2{\left( \begin{array}{c} #1 \\ #2 \end{array}\right)}

\def\vn{V_1^{\otimes n}}
\def\vnQ{V_{1,\Q}^{\otimes n}}
\def\e{\epsilon}

\def\ch{\mathfrak h^*}

\def\OO#1{\mathcal O_{#1}}
\def\Og#1{\mathcal O(#1)}
\def\Ogc#1#2{\mathcal O_{#1}(#2)}
\def\Ogp#1#2{\mathcal O(#1,#2)}
\def\Ogcp#1#2#3{\mathcal O_{#1}(#2,#3)}

\def\Lo{L_n}
\def\Ld{L_n^{\ast}}

\def\Hom{\textrm{Hom}}
\def\mc{\mathcal} 
\def\mf{\mathfrak}

\def\gvm{$\mf{p}$-Verma module} 

\def\yesnocases#1#2#3#4{\left\{
\begin{array}{ll} #1 & \mbox{#2}\\ #3 & \mbox{#4}
\end{array} \right. }

\def\drawing#1{
\begin{center} \epsfig{file=#1} \end{center} }

\def\Ti{T_i}
\def\Tui{T^i} 

\def\Mp{\mc{M}(\mf{p})}
\def\Ind{\textrm{Ind}}
\def\Cfun{\mc{C}}

\section{Introduction}

One of the most important developments in the theory of quantum groups has 
been  the discovery of canonical bases which have remarkable integrality and 
positivity properties ([Lu]). Using these bases one expects to formulate 
the representation theory of quantum groups entirely over natural numbers. 
In the case of the simplest quantum group $\U,$ such a formulation can be 
achieved via the Penrose-Kauffman graphical calculus (see [FK]). 
The positive integral structure of representation 
theory suggests that it is itself a Grothendieck ring of a certain tensor 
2-category. A strong support of this anticipation comes from 
identifying  coefficients of the transition matrix between the canonical 
and elementary bases in the $n$-th tensor power of the 
two-dimensional fundamental representation $V_1$ of $\U$ with the 
Kazhdan-Lusztig polynomials associated to $\gl$ for the maximal 
parabolic subalgebras ([FKK]). 

In this paper we will take one more step towards constructing a 
tensor 2-category with the Grothendieck ring isomorphic to the
representation category for $\U.$ The construction of tensor categories or 
2-categories with given Grothendieck groups will be referred to as 
``categorification''. We obtain a categorification of the $U(\mf{sl}_2)$ 
action in $\vn$ and the action of its commutant, the Temperley-Lieb algebra, 
using projective and Zuckerman functors between certain representation 
categories of $\gl.$ 
We extend this categorification to the comultiplication of $U(\mf{sl}_2).$ 
Our results are strongly motivated by the papers [BLM],[GrL] and 
[Gr], where the authors use the geometric rather than algebraic approach. 
In the geometric setting the categorification can be obtained via the 
categories of perverse sheaves. It is  expected that the algebraic 
and geometric languages will be equivalent, however, at the present 
moment the dictionary is still incomplete and the majority of our 
results do not allow a direct translation into the geometric 
language. Such a translation would require a nice geometric 
realization, so far unknown, 
 of singular blocks of the  highest weight categories
for $\gl$ and projective functors between these blocks. 

The main results of the paper are contained in Sections 3 and 4 and 
provide two categorifications of $U(\mf{sl}_2)$ and Temperley-Lieb algebra 
actions. Preliminary facts and 
definitions are collected in Section~\ref{sect-sl2}. 
The basis constituents of 
our construction are singular and parabolic categories of highest weight 
modules together with projective and Zuckerman functors acting on 
these categories. Projective functors in categories of $\gl$ modules, 
defined as direct summands of functors of tensoring with a 
finite-dimensional $\gl$-modules [J],[Zu], are extensively used in 
representation theory, (see  [BG] and [KV]). Being exact, projective 
functors induce linear maps in Grothendieck groups of categories 
of representations. Zuckerman functors are defined for any parabolic 
subalgebra $\mf{p}$ of $\gl$ by taking the maximal $U(\mf{p})$-locally 
finite submodule [KV]. Derived functors of Zuckerman functors are 
exact and also descend to Grothendieck groups. 
An important property of Zuckerman functors, namely their commutativity 
with the projective functors, yields in both categorifications what we 
consider the Schur-Weyl duality for $\U$ and the Temperley-Lieb algebra 
actions. 

In Section~\ref{sect-sing} we construct a categorification via 
singular blocks of the category $\O(\gl)$ of highest weight $\gl$-modules. 
More specifically, we realize $V_1^{\o n}$ as a Grothendieck group of 
the category 
\begin{equation} 
\O_n= \oplus_{k=0}^n \O_{k,n-k}, 
\end{equation} 
where $\O_{k,n-k}$ is a singular block of $\O(\gl)$ corresponding 
to the subgroup $\mathbb S_k\times \mathbb S_{n-k}$ of $\mathbb S_n.$ 
The simplest projective functors constructed by means of tensoring 
with the fundamental representation of $\gl$ and its dual descend on the 
Grothendieck group level to the action of generators $E$ and $F$ 
of $\mf{sl}_2$ (Section~\ref{functors-EF}). Various equalities between 
products of $E$ and $F$ result from functor isomorphisms 
(Section~\ref{sect-divp}). Moreover, we show that indecomposable 
projective functors in $\O_n$ correspond to elements of Lusztig 
canonical basis in the modified universal enveloping algebra 
$\dU$ (Section~\ref{canbas-indecomp}). Construction of the 
comultiplication for 
$U(\mf{sl}_2)$ requires studying the relation between categories 
$\O(\gl)\times \O(\mf{gl}_m)$ and $\O(\mf{gl}_{n+m}),$ given by 
the induction functor from the maximal parabolic subalgebra of 
$\mf{gl}_{n+m}$ that contains $\gl\oplus \mf{gl}_m.$ In particular, 
a categorification of the comultiplication formulas 
$\Delta E= E\o 1+ 1\o E$ and $\Delta F= F\o 1+ 1\o F$ for generators $E$ 
and $F$ is expressed by short exact sequences that employ certain properties 
of the induction functor (Section~\ref{sect-comult}). 

To categorify the action of the Temperley-Lieb algebra on $V_1^{\o n},$ 
we use derived functors of Zuckerman functors. 
We verify that defining relations 
for the Temperley-Lieb algebra result from appropriate functor 
isomorphisms. Projective and Zuckerman functors 
commute and that can be considered a ``functor'' version of 
the commutativity between the 
action of $\dU$ and the Temperley-Lieb algebra (Section~\ref{TL-and-Z}). 

In Section~\ref{parab-section} we construct another categorification,
 this time using parabolic subcategories of 
$\gl$ to realize $\vn$ as a Grothendieck group. There we consider the 
category 
\begin{equation} 
\O^n = \oplus_{k=0}^n \O^{k,n-k}, 
\end{equation} 
where $\O^{k,n-k}$ is a parabolic subcategory of a regular block, 
corresponding to the parabolic subalgebra of $\gl$ that contains 
$\mf{gl}_k\oplus \mf{gl}_{n-k}.$ In this picture the role of 
projective and Zuckerman functors is reversed, namely, the 
categorification of the Temperley-Lieb algebra action is given by 
projective functors, while the action of $U(\mf{sl}_2)$ is achieved via 
Zuckerman functors. We show that the composition of translation 
functors on and off the $i$-th wall at the Grothendieck group level 
yields the $i$-th generator of the Temperley-Lieb algebra by 
verifying that equivalences between these projective functors 
correspond to relations of the Temperley-Lieb algebra 
(Section~\ref{parab-tl}). This realization of the Temperley-Lieb algebra 
by functors was inspired and can be derived from the work [ES] of 
Enright and Shelton. 

In the second picture the action of $U(\mf{sl}_2)$ is 
categorified by Zuckerman functors (Section~\ref{sl2-and-Z}). 
This result can be extracted from the geometric approach of 
[BLM] and [Gr], which uses correspondences between flag varieties. 
The latter correspondences define functors between derived categories 
of sheaves, which are equivalent to the derived category of $\O^n.$ 

To summarize, we have two categorifications of the Temperley-Lieb algebra 
action on the $n$-th tensor power of the fundamental representation 
$\vn$: one by Zuckerman functors acting in singular blocks and the 
other by projective functors acting in parabolic categories. We also have 
two categorifications of the $U(\mf{sl}_2)$ algebra action on the same space: 
by projective functors between singular categories and by Zuckerman 
functors between parabolic categories. We conjecture that the 
Koszul duality functor of [BGS] exchanges these pairs of 
categorifications and that, more generally, the Koszul duality 
functor exchanges projective and Zuckerman functors. 

The categorification of the representation theory of $U(\mf{sl}_2)$ 
presented in our work explains the nature of integrality and 
positivity properties established in [FK] by a direct approach  
based on the Penrose-Kauffman calculus. 
 However, in this paper we did not reconsider 
some of the positivity and integrality results of [FK], e.g., positivity 
and integrality of the $6j$-symbol factorization coefficients. We expect 
that an extension of our approach will allow us to interpret these 
coefficients as dimensions of vector spaces of equivalences between 
appropriate functors. The problem of passing from 
 categorifying $U(\mf{sl}_2)$-representations to  $\U$-representations 
can most likely be solved by working with mixed versions of
projective and Zuckerman functors and  the 
category $\O(\gl).$ 

Moreover, many of our constructions 
admit a straightforward generalization from $U(\mf{sl}_2)$ to 
$U(\mf{sl}_m).$ In this case one should consider singular and parabolic 
categories corresponding to the subgroups $\mathbb{S}_{i_1}\times \dots 
\times \mathbb{S}_{i_m}, i_1+\dots + i_m=n,$ of $\mathbb{S}_n.$
The Temperley-Lieb algebra will be replaced by appropriate quotients 
of the Hecke algebra of $\mathbb{S}_n.$ A more difficult problem 
is to categorify the representation theory of $U(\mf{g})$ for an 
arbitrary simple Lie algebra $\mf{g}.$ Another interesting generalization 
of our results would be a categorification of the affine version of 
the Schur-Weyl duality. It is expected that in this case one should 
consider certain singular and parabolic categories of highest weight 
modules for affine Lie algebras $\widehat{\mf{gl}}_n.$
The functors of tensoring with a finite-dimensional $\mf{gl}_n$-module 
should be replaced by the Kazhdan-Lusztig tensoring with a tilting 
$\widehat{\mf{gl}}_n$-module (see [FM]). 

Finally we would like to discuss applications of our categorification 
results to a construction of topological invariants, which was the 
initial motivation for this work (see [CF]). It is well-known that 
the graphical calculus for representation theory of $\U$ and 
in particular for representations of the Temperley-Lieb algebra in 
tensor powers of $V_1$ is intimately related ([Ka]) to the Jones 
polynomial ([Jo]), which is a quantum invariant of links and can be 
extended to give an invariant of tangles. An arbitrary tangle in 
the three-dimensional space is a composition of elementary pieces such 
as braiding and local maximum and minimum tangles. To construct 
the Jones polynomial one attaches to these elementary tangles operators 
from $V_1^{\o m}$ to $\vn$ for suitable $m$ and $n$ and obtains an
isotopy invariant. 

Extending both categorifications of the Temperley-Lieb 
algebra at the end of Sections~\ref{sect-comult} and \ref{parab-tl} 
we define functors from derived categories of $\O_m$ to $\O_n$ and 
$\O^m$ to $\O^n$ corresponding to elementary tangles. Given a plane 
partition of a tangle, we can associate to it a functor, which is 
a composition of these basic functors. We conjecture that different 
plane projections produce isomorphic functors and we would get functor 
invariants of links and tangles. For links these invariants will take the 
form of $\Z$-graded homology groups. Given a diagram of a cobordism 
between two tangles, we can associate to it a natural transformation of 
functors. We expect that these natural transformations are isotopy invariants 
of tangle cobordisms, and, in the special case of a cobordism between 
empty tangles, invariants of 2-knots. To prove this conjecture 
one needs to present an arbitrary cobordism as a composition of 
elementary ones and verify all the relations between them. A complete
set of generators and relations has been found in [CS],[CRS] and was 
interpreted as a tensor 2-category in [Fi]. The match between 
tensor 2-categories arising from topology and representation theory will 
yield a graphical calculus for a categorification of the representation 
theory of $\U$ based on two-dimensional surfaces and as a consequence 
new topological invariants. 

{\it Acknowledgements: } J.B. and M.K. would like to thank the 
Institute for Advanced Study  and  the Institut des Hautes 
Etudes Scientifiques for their hospitality during 
our work on this paper. M.K. is indebted to Victor Protsak for 
interesting discussions. M.K. gratefully acknowledges NSF support in 
the form of grant DMS 9304580. The research of I.F.  was 
partially supported by NSF grant DMS 9700765.

\section{Lie algebra $\slt$ and categories of highest weight modules}
\label{sect-sl2} 

\subsection{$\dU$ and its representations}

\subsubsection{ Algebra $\dU$  }

The universal enveloping algebra of the Lie algebra $\slt$ is given by 
generators $E,F,H$ and defining relations 

$$ EF-FE =  H, \quad  HE-EH= 2E ,\quad HF-FH=-2F. $$

We will denote this algebra by $\UL.$ Throughout the paper we consider 
it as an algebra over the ring of integers $\Z.$ We will also need
two other  versions of this algebra, $\UZ$ and $\dU.$ 

Let $\UZ$ be the integral lattice in $\UL\o \Q$ spanned by 

$$E^{(a)}\binom{H}{b}F^{(c)}$$

for $a,b,c\ge 0.$ Here 

\begin{equation}
\label{divided-powers}
E^{(a)}=\frac{E^a}{a!},\quad F^{(c)}=\frac{F^c}{c!},\quad 
\binom{H}{b} = \frac{H(H-1)\dots(H-b+1)}{b!} 
\end{equation}
$E^{(a)}, F^{(a)}$ are known as \emph{divided powers} of $E$ and $F.$ 
The lattice $\UZ$ is closed under multiplication, and therefore inherits the 
algebra structure from that of $\UL\o \Q.$ Thus, $\UZ$ is an algebra over 
$\Z$ with  multiplicative generators 
$$1,E^{(a)},F^{(a)}, \binom{H}{a}, \quad a>0.$$ 

Some of the relations between the generators are written below 

\begin{eqnarray}
\label{divided-sl2-relations}
E^{(a)}E^{(b)} & = & 
\left( \begin{array}{c} a+b \\ a \end{array}\right)
E^{(a+b)} \\
F^{(a)}F^{(b)} & = & 
\left( \begin{array}{c} a+b \\ a \end{array}\right)
F^{(a+b)} \\
E^{(a)}F^{(b)} & = & 
\sum_{j=0}^{min(a,b)} F^{(b-j)}
\left( \begin{array}{c} H-a-b+2j \\ j \end{array}\right)
E^{(a-j)}. 
\end{eqnarray} 

It is easy to see that the comultiplication in $\UL\o \Q$ preserves 
the lattice $\UZ,$ i.e. $\Delta \UZ \subset \UZ \o \UZ,$ 
and the algebra $\UZ$ is actually a Hopf algebra. Note that on 
$E^{(a)},F^{(a)}$ the comultiplication is given by 

\begin{eqnarray}
\label{comult-div-E}
\Delta E^{(a)} & = & \sum_{b=0}^a E^{(b)}\o E^{(a-b)} \\
\label{comult-div-F}
\Delta F^{(a)} & = & \sum_{b=0}^a F^{(b)}\o F^{(a-b)}.  
\end{eqnarray}

Algebra $\dU$ is obtained by adjoining a system of projectors, 
one for each element of the weight lattice, to the algebra $\UZ.$ 
Start out with a $\UZ$-bimodule, freely generated by 
the set \mbox{$1_n, n\in \Z.$} Quotient it out by relations

 \begin{eqnarray} 
 \binom{H}{a}1_i & = & \binom{i}{a}1_i \\
 1_i\binom{H}{a} & = & \binom{i}{a}1_i \\
 E^{(a)}1_i & = & 1_{i+2a} E^{(a)} \\
 F^{(a)}1_i & = & 1_{i-2a} F^{(a)}.  
\end{eqnarray} 

The quotient $\UZ$-bimodule has a unique algebra structure, compatible with 
the $\UZ$-bimodule structure and such that 
\begin{equation}
1_n1_m=\delta_{n,m}1_n. 
\end{equation}

Denote the resulting $\Z$-algebra by $\dU.$ 
As a $\Z$-vector space, it is spanned by elements 
$$E^{(a)}1_n F^{(b)} \quad \mbox{ for }a,b\ge 0, n\in \Z.$$
As a left $\UZ$-module, $\dU$ decomposes into a direct sum
$$\dU =\oplus_{i\in \Z} \dU_i$$
where
$$\dU_i=\l x\in \dU|x1_i=x\r. $$ 

$\dU_i$ is spanned by $E^{(a)}1_{i-2b}F^{(b)}, a,b\ge 0.$

We will be using Lusztig's basis $\dB$ of $\dU,$ given by 
\begin{eqnarray} 
\label{canbasisexpl}
E^{(a)}1_{-i}F^{(b)} & \textrm{ for } & a,b,i\in \N, i\ge a+b \\
F^{(b)}1_i E^{(a)}   & \textrm{ for } & a,b,i\in \N, i > a+b . 
\end{eqnarray}
 \emph{Remark: }$ E^{(a)}1_{-a-b} F^{(b)} = F^{(b)}1_{a+b}
E^{(a)}.$

Define $\dB_i=\dB\cap\dU_i, i\in \Z.$ 

The important feature of Lusztig's basis is the positivity of the 
multiplication: for any  $x,y\in \dB$

$$xy=\sum_{z\in \dB}m_{x,y}^z z$$
with all structure constants $m_{x,y}^z$ being nonnegative integers. 
Although this positivity property 
 is very easy to verify, its generalization to 
quantum groups $U_q(\mathfrak g), \mathfrak g$ symmetrizable, conjectured 
by Lusztig, is, apparently, still unproved.  
The algebra $\dU$ is a special case ($q=1,\mathfrak g=\slt$) of the 
Lusztig's algebra $\dot{U}_q(\mathfrak g)$ (see [Lu]), obtained from  
the quantum group $U_q(\mathfrak g)$ by adding a system of projectors, one
for each element of the weight lattice. Lusztig defined a basis in 
$\dot{U}_q(\mathfrak g)$ and conjectured that the multiplication and 
comultiplication constants in this basis lie in $\N[q,q^{-1}].$ 
A proof of this conjecture, most likely, 
will require interpreting Lusztig's basis in terms of perverse sheaves on 
suitable varieties. 

Although in this paper we do not venture beyond $\slt,$ our results 
suggest a close link between the  Lusztig's basis of $\dot{U}_q(\mathfrak 
{sl}_N)$  and indecomposable projective functors for $\sln,$ 
$N$ and $n$ being independent parameters. 

\subsubsection{ Representations of $\dU$ } 
\label{JW-def}
Let $V_1$ be the two-dimensional representation over $\Z$ of $\UL$ spanned by 
$v_1$ and $v_{0}$ with the action of generators of $\UL$ given by 
\begin{eqnarray}
Hv_1=v_1 & Ev_1=0 & Fv_1=v_{0} \\
Hv_{0} = -v_{0} & Ev_{0}=v_1 & Fv_{0}=0 
\end{eqnarray}

In the obvious way, $V_1$ is also a representation of both $\UZ$ and $\dU.$ 
Using comultiplication, the tensor powers of $V_1$ become 
representations of $\UL, \UZ$ and $\dU.$

Denote by $V_0$ the one-dimensional representation of $\UL$ given 
by the augmentation homomorphism $\UL\to \Z$ of the universal enveloping 
algebra.  
Again, $V_0$ is a $\UZ$ and $\dU$ module in a natural way. 

Let $\delta$ be the module homomorphism $V_0\to V_1\o V_1$ given by 
\begin{equation}
\delta(1)=v_1\o v_{0}- v_{0}\o v_1   \label{eq:delta}
\end{equation}

For a sequence $I=a_1\dots a_n$ of ones and zeros, let 
$I_+$ be the number of ones in the sequence. We will denote the 
vector $v_{a_1}\o \dots \o v_{a_n}\in \vn$ by $v(I).$

We define a $\Q$-linear map (the symmetrization map) 
$p_n: \vn\o_{\Z} \Q \to \vn \o_{\Z} \Q$
by 
\begin{equation}
p_n(v(I))=
\binom{n}{I_+}^{-1}\sum_{J,J_+=I_+} v(J) \label{eq:JWprojector}
\end{equation}
where the sum on the right hand side is over all sequences 
$J$ of length $n$ with $J_+=I_+.$  
Then $p_n$ is a $\UL\o_{\Z} \Q$-module homomorphism;  in fact, it is 
the projection onto the unique $(n+1)$-dimensional irreducible 
$\UL\o_{\Z} \Q$ subrepresentation of $\vn\o_{\Z} \Q.$ 

\subsection{Temperley-Lieb algebra} 
\label{TL-algebra}

\begin{definition} 
{\rm The Temperley-Lieb algebra $TL_{n,q}$ is an algebra
over the ring $R=\Z[q,q^{-1}],$ where $q$ is a formal variable,  
with generators $U_1,\dots ,U_{n-1}$ and defining relations } 
\begin{eqnarray}
& U_iU_{i\pm 1}U_i &  = U_i \\
& U_iU_j & =U_jU_i \hspace{0.5in}|i-j|>1 \\
& U_i^2 & =-(q+q^{-1})U_i 
\end{eqnarray}
\end{definition}

The Templerley-Lieb algebra admits a geometric interpretation via
systems of arcs on the plane. Namely, as a free $R$-module, it has 
a basis enumerated by isotopy classes of systems of simple, pairwise 
disjoint arcs that connect $n$ points on the bottom of a horizontal 
plane strip with $n$ points on the top. We only consider systems
without closed arcs. Two diagrams are multiplied by concatenating 
them. If simple closed loops appear as a result
of concatenation, we remove them, each time 
multiplying the diagram by $-q-q^{-1}.$

Example: let diagrams $A$ and $B$ be as depicted below. 

\begin{center} \epsfig{file=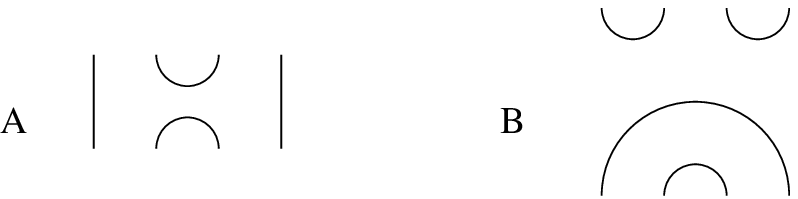} \end{center}

Then their composition $BA$ can be depicted by

\begin{center} \epsfig{file=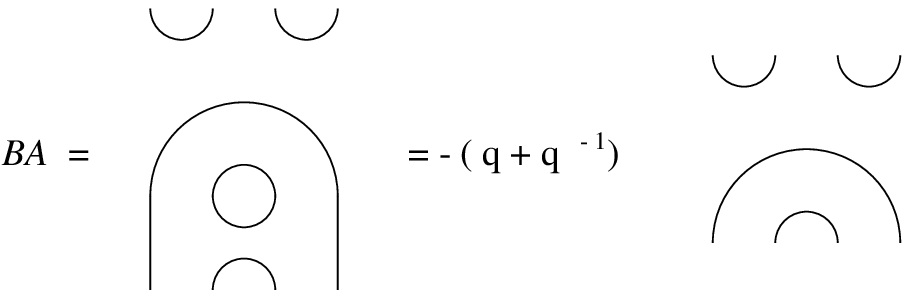} \end{center}

The generator $U_i$ of $TL_{n,q}$ is given by the diagram  

\begin{center} \epsfig{file=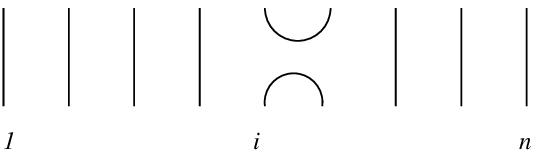} \end{center}

The defining relations have geometric interpretations. 
For instance, the 
first relation says that the diagrams below are isotopic 

\begin{center} \epsfig{file=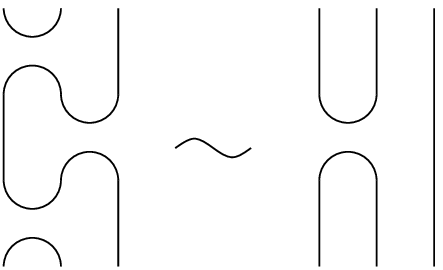} \end{center}

\begin{definition} 
{\rm The Temperley-Lieb algebra $TL_{n,1},$ 
respectively $TL_{n,-1},$ is an algebra
over the ring of integers, obtained from $TL_{n,q}$ by setting $q$ to
$1,$ respectively to $-1$,  
everywhere in the definition of the latter. } 
\end{definition} 

\vspace{0.15in} 

Thus, in $TL_{n,1}$ the value of a closed loop is $-2,$
in $TL_{n,-1}$ the value of a closed loop is $2,$  while in $TL_{n,q}$ 
a closed loop evaluates to $-q-q^{-1}.$ 

Recall that we denoted by $V_1$ the fundamental representation of 
$\UZ$. Let $u$ be an intertwiner $V_1^{\o 2}\to V_1^{\o 2}$ given 
by 
\begin{eqnarray*}
& u(v_1\o v_0)=-u(v_0\o v_1)& =v_0\o v_1-v_1\o v_0  \\
& u(v_1\o v_1)=u(v_0\o v_0) & = 0 
\end{eqnarray*}
Then $V_1^{\o n}$ is a representation of $TL_{n,1}$ with $U_i$ 
acting by $Id^{\o (i-1)}\o u \o Id^{\o (n-i-1)}.$ 
This action commutes with the Lie algebra $\slt$ action on the same 
space. 

The Temperley-Lieb algebra allows a generalization into the so-called 
 Temperley-Lieb category, as we now explain (for more details, see 
 [KaL],[Tu]).

\begin{definition} {\rm The Temperley-Lieb category 
$TL$ has objects enumerated by nonnegative integers: $\mbox{Ob}(TL)=
\l \overline{0},\overline{1},\overline{2},\dots \r .$ 
The set of morphisms from $\overline{n}$ to $\overline{m}$ 
is a free $R$-module with a basis over $R$ given by the isotopy classes 
of systems of $\frac{n+m}{2}$ 
simple, pairwise disjoint arcs inside a horizontal 
strip on the plane that connect in pairs $n$ points on the bottom and $m$ 
points on the top in some order. 

Morphisms are composed by concatenating their diagrams. If closed loops 
appear after concatenation, we remove them, multiplying the diagram 
by $-q-q^{-1}$ to the power equal to the number of closed loops. } 
\end{definition} 

\vspace{0.15in} 

An example of a morphism from $\overline{5}$ to $\overline{3}$ is 
depicted below. 

\begin{center} \epsfig{file=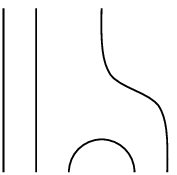} \end{center}
If $n+m$ is odd, there are no morphisms from 
$\overline n$ to $\overline m.$ 
Denote by $\cap_{i,n}$ for  $n\ge 2, 1\le i\le n-1$ the morphism of $TL$
from $\overline{n}$ to $\overline{n-2}$ given by the following 
diagram

\drawing{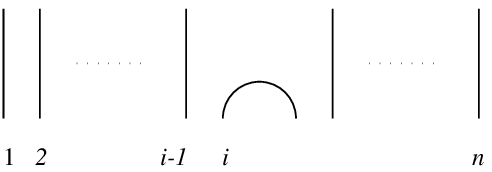}

The diagram consists of $n-1$ arcs.  One of the arcs connects the 
$i$-th bottom point (counting from the left) with the $(i+1)$-th bottom 
point. The remaining arcs connect the $k$-th bottom point for $1\le k<i$ with 
the $k$-th top point and the $k$-th bottom point for $i+2\le k\le n$ with 
the $(k-2)$-th top point. 

Denote by $\cup_{i,n}, n\ge 0, 1\le i\le n+1$ the morphism 
in $TL$ from $\overline{n}$ to $\overline{n+2}$ given by the diagram

\drawing{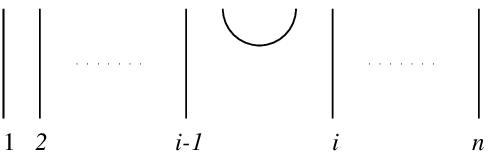}

Denote by $\mbox{Id}_{\overline{n}}$ the identity morphism from
$\overline{n}$ to $\overline{n}.$ This morphism  can be depicted by a diagram 
that is made of $n$ vertical lines: 

\drawing{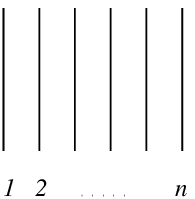}

The morphisms $\cap_{i,n}$ and $\cup_{i,n}$ will serve as generators 
of the set of morphisms in the Temperley-Lieb category. The following 
is a set of defining relations for $TL$  

\begin{eqnarray}
& \cap_{i+1,n+2}\circ \cup_{i,n} & = \mbox{Id}_{\overline{n}} 
\label{first-relation}\\
& \cap_{i,n+2}\circ \cup_{i+1,n} & = \mbox{Id}_{\overline{n}} 
\label{second-relation} \\
& \cap_{j,n}\circ \cap_{i,n+2}  &  = \cap_{i,n}\circ \cap_{j+2,n+2} 
\quad i\le j \label{third-relation} \\
& \cup_{j,n-2}\circ \cap_{i,n}  &  = \cap_{i,n+2}\circ \cup_{j+2,n} 
\quad i\le j \\
& \cup_{i,n-2}\circ \cap_{j,n}  &  = \cap_{j+2,n+2}\circ \cup_{i,n} 
\quad i\le j \\
& \cup_{i,n+2}\circ \cup_{j,n}  &  = \cup_{j+2,n+2}\circ \cup_{i,n} 
\quad i\le j 
\label{sixth-relation}\\
&  \cap_{i,n+2}\circ \cup_{i,n} & = - (q+q^{-1})\mbox{Id}_{\overline{n}}
\label{last-relation}
\end{eqnarray} 

The first $6$ types of relations come from isotopies of certain pairs 
of diagrams. For example, relations (\ref{first-relation}) and 
(\ref{third-relation}) correspond to the isotopies 

\begin{center} \epsfig{file=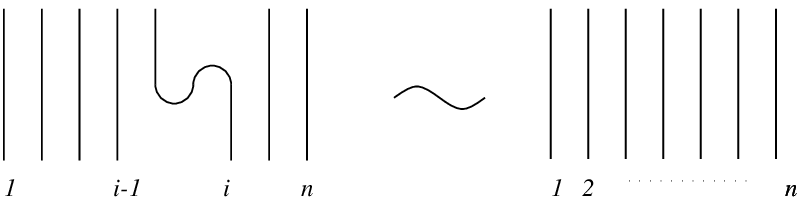} \hspace{.5in} \end{center}
and 
\begin{center} \epsfig{file=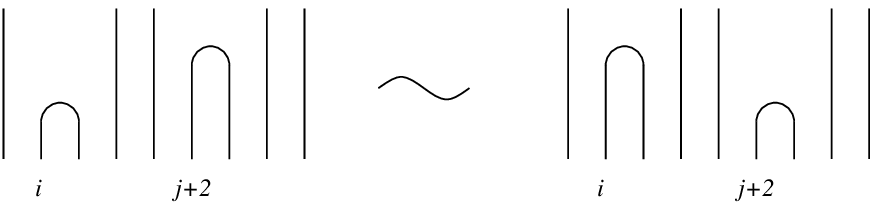} \hspace{.5in} \end{center}
respectively. 
Algebra $TL_{n,q}$ is the algebra of endomorphisms of the 
object  $\overline n$  of the Temperley-Lieb category. 

\subsection{The category of highest weight modules over a  reductive 
Lie algebra} 

\subsubsection{ Definitions } 
\label{definitions-LA}

In this section all Lie algebras and their representations are 
defined over the field $\C$ of complex numbers. 
Let $\mf{g}$ be a finite-dimensional reductive Lie algebra and 
$U(\mf{g})$ its universal enveloping algebra.  
Fix a 
triangular decomposition $\mf{g}=\mf{n}_+\oplus \mf{h}\oplus \mf{n}_-.$ 
Let $R_+$ be the set of positive roots and 
$\rho$ the half-sum of positive roots. 
For $\lambda\in \ch$ denote by $M_{\lambda}$ the Verma module with
highest weight $\lambda-\rho$ and by $L_{\lambda}$ the irreducible quotient
of $M_{\lambda}.$
The module $L_{\lambda}$ is finite-dimensional if and only if $\lambda-\rho$ 
is an integral dominant weight. 

Denote by $\Og{\mf{g}}$ the category of finitely generated $U(\mf{g})$-modules 
that are $\mf{h}$-diagonalizable and locally $U(\mf{n}_+)$-nilpotent. 
The category $\Og{\mf{g}}$ is called \textit{ the category of 
highest weight $\mf{g}$-modules}. 
Let $P_{\lambda}$ denote the projective 
cover of $L_{\lambda}$ (see [BGG] for the existence of projective covers). 

If $A$ is an additive category, denote by $K(A)$ the Grothendieck group 
of $A.$ Denote by $[M]$ the image of an object $M\in \textrm{Ob}(A)$ 
in the Grothendieck group of $A.$ Denote by $D^b(A)$ the bounded derived 
category of an abelian category $A.$

\subsubsection{ Projective functors }

This section is a brief introduction to projective 
functors. We refer the reader to Bernstein-Gelfand paper [BG] for a 
detailed treatment and further references. 

Denote by $\Theta$ the set of maximal ideals of the center $Z$ of 
$U(\mf{g}).$ 
We can naturally identify $\Theta$ with the quotient of the weight space
$\ch$ by the action by reflections 
of the Weyl group $W$ of $\gl.$ 
We will denote by $\eta$ the quotient map $\ch\to \Theta.$ 
For $\theta\in\Theta$ denote by $J_{\theta}$ the corresponding maximal 
ideal of $Z.$ Thus, the Verma module $M_{\lambda}$ with the highest weight 
$\lambda - \rho$ is annihilated by the maximal central 
ideal $J_{\eta(\lambda)}.$ 

For $\theta\in \Theta$ denote by $\Ogc{\theta}{\mf{g}}$ a full subcategory of 
$\Og{\mf{g}}$ consisting of modules that are annihilated by some power of 
the central ideal $J_{\theta}:$ 

\begin{equation} 
M\in \Ogc{\theta}{\mf{g}} \hspace{0.2in} \Longleftrightarrow 
 \hspace{0.2in} M\in \Og{\mf{g}} \mbox{ and }J_{\theta}^N M =0 
\textrm{ for sufficiently large }N
\end{equation}
A module $M\in \O(\mf{g})$ belongs to $\Ogc{\theta}{\mf{g}}$ if and only 
if all of the simple subquotients of $M$ are isomorphic to simple modules 
$L_{\lambda}, \lambda \in \eta^{-1}(\theta).$ 
We will call modules in $\Ogc{\theta}{\mf{g}}$ 
\textit{highest weight modules 
with the generalized central character} $\theta$. 
The category $\Og{\mf{g}}$ splits as a direct sum of categories 
$\Ogc{\theta}{\mf{g}}$ 
over all $\theta\in \Theta.$ 

Denote by $\textrm{proj}_{\theta}$ the functor from $\Og{\mf{g}}$ to 
$\Ogc{\theta}{\mf{g}}$ 
that, to a module $M,$ associates the largest submodule of $M$
with the generalized central character $\theta.$ 
Let $F_V$ be the functor of tensoring with a finite-dimensional 
$\mf{g}$-module $V.$ 

\begin{definition} {\rm $F:\Og{\mf{g}} \to \Og{\mf{g}}$ is a projective 
functor if it is isomorphic to a direct summand of the functor $F_V$ for some 
finite dimensional module $V.$ } 
\end{definition}

\vspace{0.15in} 

The functor $\textrm{proj}_{\theta}$ is an example of a 
 projective functor, since it is  
a direct summand of the functor of tensoring with the  
one-dimensional representation. We have an isomorphism of functors 
\begin{equation} 
F_V=\oplusop{\theta_1,\theta_2\in \Theta}
(\textrm{proj}_{\theta_1}\circ F_V\circ \textrm{proj}_{\theta_2})
\end{equation}

Any projective functor takes projective objects 
 in $\Og{\mf{g}}$ to projective objects.  
The composition of projective functors is again a projective functor. 
Each projective functor splits as a direct sum of 
indecomposable projective functors. 

Projective functors are exact. 
Therefore, they induce endomorphisms of the Grothendieck 
group of the category $\Og{\mf{g}}.$ The following result is proved in [BG]: 

\begin{prop}
\label{prop:projective-functors-isomorphic}
Let $\lambda$ be a dominant integral weight, $\theta=\eta(\lambda)$ 
 and 
$F,G$ projective functors from $\Ogc{\theta}{\mf{g}}$ to 
$\Og{\mf{g}}.$ Then 

\begin{enumerate}
\item Functors 
$F$ and $G$ are isomorphic if and only if 
the endomorphisms of $K(\O(\mf{g}))$ induced by 
$F$ and $G$ are equal. 
\item Functors $F$ and $G$ are isomorphic if and only if 
 modules $FM_{\lambda}$ and $GM_{\lambda}$ are isomorphic.  
\end{enumerate}
\end{prop}

We will be computing the action of projective functors on Grothendieck 
groups of certain subcategories of the category of highest weight 
modules. The simplest basis in the Grothendieck group of $\O(\mf{g})$ 
is given by images of Verma modules. The following 
proposition shows that this basis is also handy for 
writing the action of projective functors on the Grothendieck group
of $\O(\mf{g}).$ 

\begin{prop}
\label{prop:Verma-tensor-finite}
Let $V$ be a finite-dimensional $\mf{g}$-module,
$\mu_1,\dots,\mu_m$ a multiset of weights of $V,$
 $M_{\chi}$ the Verma
module with the highest weight $\chi-\rho,$ then 

\begin{enumerate}
\item The module $V\o M_{\chi}$ admits a fitration with 
successive quotients isomorphic to Verma modules 
$M_{\chi+\mu_1},\dots ,M_{\chi+\mu_m}$ (in some order).
\item We have an equality in the Grothendieck group 
$K(\Og{\mf{g}})$:
$$[V\otimes M_{\chi}]=\sum_{i=1}^m [M_{\chi+\mu_i}]$$
\end{enumerate}
\end{prop}

\subsubsection{ Parabolic categories }
\label{parabolic-categories} 

Let $\mf{p}$ be a parabolic subalgebra of $\mf{g}$ that contains 
$\mf{n}_+\oplus \mf{h}.$ Denote by $\Ogp{\mf{g}}{\mf{p}}$ the full 
subcategory of $\Og{\mf{g}}$ that consists of $U(\mf{p})$ locally 
finite modules. Notice that projective functors preserve 
subcategories  $\Ogp{\mf{g}}{\mf{p}}.$

A generalized Verma module relative to a parabolic subalgebra 
$\mf{p}$ of $\mf{g}$ (see [Lp],[RC]) will be called 
a $\mf{p}$-Verma module. The Grothendieck group of $\Ogp{\mf{g}}{\mf{p}}$ 
is generated by images $[M]$ of generalized Verma modules.

For a central character $\theta$ and a parabolic subalgebra $\mf{p}$ of 
$\mf{g},$ denote by $\O_{\theta}(\mf{g},\mf{p})$ the full subcategory of 
$\O(\mf{g})$ consisting of $U(\mf{p})$-locally finite modules 
annihilated by some power of the central ideal $J_{\theta}.$ The 
category  $\O_{\theta}(\mf{g},\mf{p})$  is the intersection of 
subcategories  $\O_{\theta}(\mf{g})$  and 
$\O(\mf{g},\mf{p})$ of $\O(\mf{g}).$

The following lemma is an obvious
generalization of a special case of Lemma 3.5 in [ES].

\begin{lemma}
\label{functor-equivalence}
 Let $T,S$ be covariant, exact functors from
$\Ogcp{\theta}{\mf{g}}{\mf{p}}$ to some abelian category $\mc A$ and let
$f$ be a natural transformation from
$T$ to $S.$ If $f_M: T(M)\to S(M)$ 
is an isomorphism for each generalized Verma
module $M\in \Ogcp{\theta}{\mf{g}}{\mf{p}},$ then $f$ is an isomorphism 
of functors.
\end{lemma}

Let $\mf{g}_1,\mf{g}_2$ be reductive Lie algebras, with fixed 
Cartan subalgebras 
$\mf{h}_j\subset \mf{g}_j,j=1,2.$ Suppose we have two parabolic 
subalgebras $\mf{p}_1,\mf{p}_2$ such that
$\mf{h}_j\subset \mf{p}_j\subset \mf{g}_j.$
Fix central characters $\theta_j$ of $\mf{g}_j.$ 

\begin{lemma} 
\label{category-equivalence}
Suppose that we have two exact functors
\begin{eqnarray*}
& f_{12}: \Ogcp{\theta_2}{\mf{g}_2}{\mf{p}_2} & \to 
\Ogcp{\theta_1}{\mf{g}_1}{\mf{p}_1} \\
& f_{21}: \Ogcp{\theta_1}{\mf{g}_1}{\mf{p}_1} & \to 
\Ogcp{\theta_2}{\mf{g}_2}{\mf{p}_2} 
\end{eqnarray*}
such that $f_{21}$ is isomorphic to both left and right adjoint 
functors of $f_{12},$ $f_{21}$ takes $\mf{p}_1$-Verma modules to 
$\mf{p}_2$-Verma modules, $f_{12}$ takes $\mf{p}_2$-Verma modules to 
$\mf{p}_1$-Verma modules, 
$f_{12}f_{21}(M)$ is isomorphic to $M$ for any $\mf{p}_1$-Verma module M and 
$f_{21}f_{12}(M)$ is isomorphic to $M$ for any $\mf{p}_2$-Verma module M. 
Then $f_{12}$ and $f_{21}$ are equivalences of categories and 
the natural transformations 
$$a:f_{12}f_{21}\to Id,\quad b: Id \to f_{12}f_{21}$$ 
coming from the adjointness are isomorphisms of functors. 
\end{lemma}

\emph{Proof: }By Lemma~\ref{functor-equivalence} it suffices to prove that, 
for any $\mf{p}_1$-Verma module $M,$ the module morphisms 
$$a_M: f_{12}f_{21}(M)\to M,\quad b_M: M\to f_{12}f_{21}(M)$$
are isomorphisms. Note that $f_{21}(M)$ is a $\mf{p}_2$-Verma 
module and $f_{12}f_{21}(M)$ is isomorphic to $M.$ 
The hom space $\mbox{Hom}_{\mf{g}_1}(M,M)$ is one-dimensional
and all morphisms are just scalings of the identity morphism. We have
a natural isomorphism
$$
\mbox{Hom}_{\mf{g}_1}(f_{12}f_{21}(M),M) = 
\mbox{Hom}_{\mf{g}_2}(f_{21}M,f_{21}M).$$
Under this isomorphism, the map $a_M: M \to M$ corresponds
to the identity map $f_{21}M\to f_{21}M.$ 
This identity map generates the space
$\mbox{Hom}_{\mf{g}_2}(f_{21}M,f_{21}M),$ therefore $a_M$ generates 
the hom space $\mbox{Hom}_{\mf{g}_1}(M,M),$ and 
thus, $a_M$ is an isomorphism of $M$,
being a non-zero multiple of the identity map. Therefore, $a$ is an 
isomorphism of functors. 
Similarly, $b$ is a functor isomorphism. $\square$

This lemma is used in Section~\ref{TL-and-Z} in the proof of 
Theorem~\ref{functors-and-isotopies}.

\subsubsection{Zuckerman functors}

Here we recall the basic properties of Zuckerman derived functors,
following [ES] and [EW]. Knapp and Vogan's book [KV] contains a 
complete treatment of Zuckerman functors, but here we will only need some 
basic facts. 

Throughout the paper we restrict Zuckerman functors to the category 
of highest weight modules. 

Let $\mf{g},\mf{p}$ be as in Section~\ref{parabolic-categories}. 
The parabolic Lie algebra 
$\mf{p}$ decomposes as a direct sum $\mf{m}\oplus \mf{u}$ where $\mf{m}$ 
is the maximal reductive subalgebra of $\mf{p}$  
 and $\mf{u}$ is the nilpotent radical of $\mf{p}.$  The reductive subalgebra 
$\mf{m}$ contains the Cartan subalgebra $\mathfrak h$ of $\mf{g}.$ 
Let $d=\mbox{dim}(\mf{m})-\mbox{dim}(\mf{\mathfrak h}).$
We denote by $\ast$ the contravariant duality functor  in $\O.$

Let $\Gamma_{\mf{p}}$ be the functor from $\Og{\mf{g}}$ to 
$\Ogp{\mf{g}}{\mf{p}}$ that to a module $M\in \Og{\mf{g}}$ associates 
its maximal locally $U(\mf{p})$-finite submodule. 
$\Gamma_{\mf{p}}$ is called the Zuckerman functor. Functor 
$\Gamma_{\mf{p}}$ is left exact and the category $\Og{\mf{g}}$ has 
enough injectives, so we can define the derived functor 
$$\mc{R}\Gamma_{\mf{p}} : D^b(\Og{\mf{g}}) \lra D^b(\Ogp{\mf{g}}{\mf{p}})$$
and its cohomology functors 
$$\mc{R}^i\Gamma_{\mf{p}}: \Og{\mf{g}} \lra \Ogp{\mf{g}}{\mf{p}}.$$

\begin{prop}
\label{Zproperties}
\begin{enumerate}
\item For $i>d , \mc{R}^i\Gamma_{\mf{p}}=0.$
\item Projective functors commute with Zuckerman functors. More precisely,
if $F$ is a projective functor, then there are natural isomorphisms
 of functors 
\begin{eqnarray*}
F\circ \Gamma_{\mf{p}} & \cong & \Gamma_{\mf{p}}\circ F \\
F\circ \mc{R}\Gamma_{\mf{p}} & \cong & \mc{R}\Gamma_{\mf{p}}\circ F
\end{eqnarray*}
\item The functors $M\longmapsto \mc{R}^i\Gamma_{\mf{p}}(M)$ and
$M\longmapsto \mc{R}^{d-i}\Gamma_{\mf{p}}
(M^{\ast})^{\ast}, M\in \O(\mf{g})$ are 
naturally equivalent.
\item $\mc{R}^d\Gamma_{\mf{p}}$ is isomorphic to the functor that to a 
module $M\in \O(\mf{g})$ associates the maximal locally 
$\mf{p}$-finite quotient of $M.$

\end{enumerate}
\end{prop}

\emph{Proof: } See [EW].  Zuckerman functors commute with functors 
of tensor product by a finite-dimensional module. A projective functor 
is a direct summand of a tensor product functor. Part 2 of the 
proposition follows. $\square$

\subsection{ Singular blocks of the highest weight category 
for $\gl$}

\subsubsection{ Notations }
\label{section-def-gln}

We fix once and for all a triangular decomposition
$\mathfrak{n}_+\oplus \mathfrak{h}\oplus \mathfrak{n}_-$ of 
the Lie algebra $\gl.$ The Weyl group 
of $\gl$ is isomorphic to the symmetric group $S_n.$  
Choose an orthonormal basis $e_1,\dots,e_n$ in the Euclidean space 
$\R^n$ and identify the complexification $\C\o_{\R}\R^n$ with 
the dual $\ch$ of Cartan subalgebra so that 
$R_+=\l e_i-e_j, i < j\r$ is  the set of positive roots and 
$\alpha_i=e_i-e_{i+1},1\le i\le n-1$ are simple roots.
The generator $s_i$ of the Weyl group  $W=S_n$ acts on $\ch$ by permuting
 $e_i$ and $e_{i+1}.$
Denote by $\rho$ the half-sum of positive roots 
$$\rho=\frac{n-1}{2}e_1+\frac{n-3}{2}e_2+\dots+\frac{1-n}{2}e_n$$
Sometimes we will use the notation $\rho_n$ instead of $\rho.$ 

For a sequence $a_1,\dots,a_n$ of zeros and ones, denote by 
$M(a_1\dots a_n)$ the Verma module with the 
highest weight $a_1e_1+\dots+a_ne_n-\rho.$ Similarly, $L(a_1\dots a_n)$
will denote the simple quotient of $M(a_1\dots a_n)$ and $P(a_1\dots a_n)$ 
the minimal projective cover of $L(a_1\dots a_n)$. 

The sequence of $n$ zeros, respectively ones, will be denoted by $0^n$, 
respectively $1^n.$ If $I_1,I_2$ are two sequences of $0$'s and $1$'s, 
we denote their concatenation by $I_1I_2.$ 

Recall that $\lambda_i=e_1+\dots+e_i$ is a fundamental 
weight of $\gl$. Denote by $\theta_i=\eta(\lambda_i)$ the 
corresponding central 
character of $\gl.$ We denote the category $\O_{\theta_i}(\gl)$ by
 $\O_{i,n-i}.$ 
A module $M\in \O(\gl)$ lies in $\O_{i,n-i}$ if and only if all of its simple 
subquotients are isomorphic to $L(a_1\dots a_n)$ for sequences 
$a_1\dots a_n$ of zeros and ones with exactly $i$ ones. 

Denote by $\O_n$ the direct sum of categories $\O_{i,n-i}$ as 
$i$ ranges over all integers from $0$ to $n$: 
\begin{equation} 
\O_n= {\mathop{\oplus}\limits_{i=0}^n}\O_{i,n-i}. 
\end{equation}   
When $i<0$ or $i> n,$ denote by $\O_{i,n-i}$ the subcategory of $\O(\gl)$ 
consisting of the zero module. 
We have an isomorphism of Grothendieck groups 
\begin{equation} 
K(\O_n)= {\mathop{\oplus}\limits_{i=0}^n}K(\O_{i,n-i}).
\end{equation} 

Let $\Upsilon_n$ be 
the isomorphism of abelian groups $\Upsilon_n: K(\O_n)\to \vn$
given by 
\begin{equation}
\label{iso-upsilon}  
\Upsilon_n [M(a_1\dots a_n)]=v_{a_1}\o \dots \o v_{a_n}.
\end{equation}

\subsubsection{ Simple and projective module bases in the Grothendieck 
group $K(\O_n)$} 
\label{sec:canbasis-formulas}

Abelian group isomorphism $\Upsilon_n$ identifies the Grothendieck 
group of $\O_n$ and the $\dU$-module $\vn.$ In this section we describe 
the images under $\Upsilon_n$ 
of simple modules and indecomposable projectives in $\O_n.$ 
The only result of this section that we use later is the 
formula (\ref{fla:0101})
for the image of the indecomposable projective $P(0^j1^k0^l1^m).$ 
This formula is 
used in the proof of Theorem~\ref{thrm:indecomposable-proj-functors}.

Let us define $3$ bases in $\vn.$ We are working over $\Z$ and thus 
$\vn$ is a free abelian group of rank $2^n.$ 
We will parametrize basis elements 
by sequences of ones and zeros of length $n.$ 

First, the basis $\l v(a_1\dots a_n),a_i\in \l 0,1\r \r$ 
will be given by 
$$v(a_1\dots a_n)=v_{a_1}\o \dots \o v_{a_n}$$
so that, for a sequence $I$ of length $n$ of zeros and ones, we 
have $\Upsilon_n[M(I)]=v(I).$ We will call this basis \textit{ the product 
basis of } $\vn.$  

Next we introduce the 
basis $\l l(a_1\dots a_n),a_i\in \l 0,1\r  \r$ 
by induction of $n$ as follows: 

(i) $l(1)=v_1,l(0)=v_{0},$ 

(ii) $l(0a_2\dots a_n)=v_{0}\o l(a_2\dots a_n),$

(iii) $l(a_1 \dots a_{n-1}1)=l(a_1,\dots,a_{n-1})\o v_1,$

(iv) 
\begin{eqnarray*}
&  & l(a_1\dots a_{i-1}10a_{i+2}\dots a_n) = \\
& & (\mbox{Id}^{\o (i-1)}\o \delta \o \mbox{Id}^{\o (n-i-1)})l(a_1\dots a_{i-1}
a_{i+2}\dots a_n),
\end{eqnarray*}
where $\delta$ is the intertwiner $V_0\to V_1\o V_1$ defined by 
the formula (\ref{eq:delta}) and $\mbox{Id}$ denotes the identity homomorphism 
of $V_1.$  
These rules are consistent and uniquely define $l(a_1\dots a_n)$ 
for all $a_1\dots a_n.$ 

Let us now define the basis 
$\l p(a_1\dots a_n),a_i\in \l 0,1\r \r.$ 
Let $I$ be a sequence of zeros and ones. Then define the basis inductively by 
the rules 

(i) $p(1I)=v_1\o 
p(I),$

(ii) $p(I0)=p(I)\o v_0,$

(iii) If a sequence $I_1$ is either empty or ends with at least 
$j$ zeros, and a sequence $I_2$ is either empty or starts with 
at least $k$ ones, then 
$$p(I_10^j1^kI_2)=\binom{j+k}{j}
(\mbox{Id}^{\o |I_1|}\o p_{j+k}\o \mbox{Id}^{\o |I_2|} ) 
p(I_11^{k}0^{j}I_2)$$
where $|I|$ stands for the length of $I.$ 

\begin{prop} These rules are consistent, and for each sequence $I$ define an 
element of $V_1^{\o |I|}.$ 
\end{prop}

Proposition follows from the results of [K], Section 3,
 setting $q$ to $1.$ 
Note that it is not even obvious that $p(I)$ lies in 
$\vn$ because the projector $p_i$ is defined (see 
Section~\ref{JW-def}) 
as an operator in  $V_1^{\o i}\o_{\Z} \Q$, rather than in $V_1^{\o i}.$

\begin{prop}
\label{prop:simple-projective-bases}
 The isomorphism $\Upsilon_n:K(\O_n)\to \vn$ takes the images of simple, 
resp. indecomposable projective modules 
to elements of the basis $\l l(I) \r_I $, resp. 
$\l p(I)\r_I$ of $\vn:$

\begin{eqnarray*}
\Upsilon_n[L(I)] & = & l(I) \\
\Upsilon_n[P(I)] & = & p(I)
\end{eqnarray*}
\end{prop}

\emph{Proof: } The rules (i)-(iii) for $p(I)$ 
can be used to write down the relations 
between the coefficients of the transformation matrix 
from the basis $\l v(I)\r_I$ to the basis $\l p(I) \r_I$ of 
$\vn.$ These relations are equivalent to the $q=1$ 
specialization of Lascoux-Sch\"{u}tzenberger's 
recursive formulas (see [LS] and [Z]) 
for the Kazhdan-Lusztig polynomials in the Grassmannian case, as follows 
from the computation at the end of [FKK] (again, setting $q$ to $1$). 

Kazhdan-Lusztig polynomials 
in the Grassmannian case for $q=1$ are 
coefficients of the transformation matrix from the Verma module to
the simple module basis of Grothendieck groups of certain parabolic 
subcategories of a regular block of $\O(\gl).$ These parabolic subcategories 
are Koszul dual (see [BGS], Theorem 3.11.1 for the general statement) 
to the singular blocks $\O_{i,n-i}, 0\le i\le n$ of 
$\O(\gl).$ Koszul duality functor descends to the isomorphism 
of Grothendieck groups 
that exchanges simple and projective modules in corresponding categories. 
Therefore, Lascoux-Sch\"{u}tzenberger's formulas also describe 
coefficients of decomposition of projective modules in the Verma module 
basis of $K(\O_n).$ It now follows that $\Upsilon_n[P(I)]=p(I).$ 

Introduce a bilinear form on $K(\O_n)$ by $<[M(I)],[M(J)]>=\delta_I^J.$ 
The BGG reciprocity implies $<[P(I)],[L(J)]> =\delta_I^J,$ i.e., 
the basis $\l [L(I)] \r $  of $K(\O_n)$ is dual to the basis 
$\l [P(I)] \r$ of $K(\O_n).$ Abelian group isomorphism $\Upsilon_n$ 
transform this bilinear form to a 
bilinear form on $\vn$  such that $<v(I),v(J)>=\delta_I^J.$ 
From the main computation of Chapter 3 of [K], specializing $q$ to $1$, 
it follows 
that bases $\l l(I) \r $ and 
$\l p(I) \r $ are orthogonal relative to this form. 
We have  $\Upsilon_n[M(I)]=v(I)$ by definition of $\Upsilon_n$ and 
we have already established that $\Upsilon_n[P(I)]= p(I)$. We conclude that 
$\Upsilon_n[L(I)]= l(I).$ 

$\square$

To prove Theorem~\ref{thrm:indecomposable-proj-functors} 
we will need an explicit formula for $p(0^j1^k0^l1^m):$

\begin{equation}
\label{fla:0101}
\begin{array}{lll}
& & p(0^j1^k0^l1^m) = \\
& & =  \binom{j+k}{k}\binom{j+l+m}{m} 
(p_{j+k}\o \mbox{Id}^{\o (l+m)})(\mbox{Id}^{\o k}\o p_{j+l+m}) \\
& & (v_1^{\o (k+m)}\o v_0^{\o (j+l)}) \textrm{ if }k\le l \\
& & = \binom{l+m}{m}\binom{j+k+m}{j}
(\mbox{Id}^{\o (j+k)}\o p_{l+m})(p_{j+k+m}\o \mbox{Id}^{\o l}) \\
& & (v_1^{\o (k+m)}\o 
v_0^{\o (j+l)}) \textrm{ if }k\ge l. 
\end{array}
\end{equation}

This formula follows from recurrent relations (i)-(iii) 
for $p(I)$ that we gave earlier in this section. 


\section{Singular categories} 
\label{sect-sing} 

\subsection{ Projective functors and 
$\mathfrak{sl}_2$ }
\label{proj-sl2}

\subsubsection {Projective functors $\E$ and $\F$} 
\label{functors-EF}
 
Let $\Lo$ be the $n$-dimensional representation of $\gl$ with 
vweights $e_1,e_2,\dots ,e_n.$ The dual representation $\Ld$ has 
weights $-e_1,-e_2,\dots, -e_n.$ 

Recall that $\O_{i,n-i}, i=0,1,\dots n$ 
is the singular block of $\O(\gl)$ consisting 
of modules with generalized central character $\theta_i=\eta(\lambda_i).$ 
For $i<0$ and $i>n$ we defined $\O_{i,n-i}$ to be the trivial subcategory 
of $\O(\gl).$ 

Denote by $\E_i$ the projective functor 
$$(\textrm{proj}_{\theta_{i+1}})
\circ F_{\Lo}: \quad \O_{i,n-i}\to \O_{i+1,n-i-1}$$
given by tensoring with the $n$-dimensional representation
$\Lo$ and then taking the 
largest submodule of this tensor product that lies in $\O_{i+1,n-i-1}.$ 

Similarly,  
denote by $\F_i$ the projective functor 
from $\O_{i.n-i}$ to $\O_{i-1,n-i+1}$
given by tensoring with $\Ld$ and then taking the largest  
submodule that belongs to $\O_{i-1,n-i+1}$

\begin{theorem}
\label{thrm:EF-relation}
For $i=0,1,\dots n$ there are isomorphisms of projective functors 
\begin{equation} 
(\E_{i-1}\circ \F_i) \oplus
\textrm{Id}^{\oplus(n-i)}
\cong
(\F_{i+1}\circ \E_i) \oplus
\textrm{Id}^{\oplus i}  
\end{equation} 
where $Id$ denotes the identity functor $Id:\O_{i,n-i}\to \O_{i,n-i}.$ 
\end{theorem}

\emph{Proof:} By Proposition~\ref{prop:projective-functors-isomorphic}
 it suffices to check the equality
of endomorphisms of the Grothendieck group of $\O_{i,n-i}$ 
induced by these projective functors.

Denote by $[\E_i]$ and $[\F_i]$ 
the homomorphisms of the Grothendieck groups induced by functors 
$\E_i$ and $\F_i:$
\begin{eqnarray*}
& [\E_i]: &  
K(\O_{i,n-i})\to K(\O_{i+1,n-i-1}) \\
& [\F_i]: & 
K(\O_{i,n-i})\to K(\O_{i-1,n-i+1})
\end{eqnarray*}
The Grothendieck group $K(\O_{i,n-i})$ is
free abelian of rank $\binom{n}{i}$ and 
is spanned by images $[M(a_1\dots a_n)]$ 
of Verma modules $M(a_1\dots a_n)$ for all possible sequences
$a_1\dots a_n$ of zeros and ones with $i$ ones.

By Proposition~\ref{prop:Verma-tensor-finite}
\begin{equation}
[M(a_1\dots a_n)\otimes \Lo]=
\sum_{j=1}^n[M(a_1\dots a_j'\dots a_n)]
\end{equation}
where $a_j'=a_j+1.$ 

The functor $\E_i$ is a composition of tensoring with $\Lo$
 and a projection onto $\O_{i+1,n-i-1},$ hence we get 

\begin{prop} Let $a_1\dots a_n$ be a sequence of zeros and ones that 
contains $i$ ones. Then 

\begin{equation}
 [\E_iM(a_1 \dots a_n)] =
\sum_{j=1,a_j=0}^{n}[M(a_1\dots a_{j-1}1 a_{j+1}\dots a_n)]
\end{equation}

\end{prop}

In the same fashion, we obtain 

\begin{prop} Let $a_1\dots a_n$ be a sequence of zeros and ones that 
contains $i$ ones. Then 

\begin{equation}
[\F_iM(a_1 \dots a_n)] =
\sum_{j=1,a_j=1}^{n}[M(a_1\dots a_{j-1}0 a_{j+1}\dots a_n]
\end{equation}

\end{prop}

Therefore, after identifying the Grothendieck group $K(\O_n)$ with 
$\vn$ via the isomorphism $\Upsilon_n$ (formula (\ref{iso-upsilon}), 
we see that 
maps $[\E_i],[\F_i]$  from $K(\O_{i,n-i})$ to $K(\O_{i+ 1, n-i-1})$
and $K(\O_{i-1,n-i+1})$    
coincide with maps induced by the 
Lie algebra $\slt$ generators $E$ and $F$ on 
the weight $2i-n$ subspace of the module $\vn.$  
This immediately gives 
\begin{prop} We 
have the following equality of endomorphisms of the abelian group 
$K(\OO{i,n-i})$ : 
$$[\E_{i-1}] [\F_i]+(n-i)\textrm{Id} = [\F_{i+1}] [\E_i]+ 
 i \cdot \textrm{Id}$$  
\end{prop} 

Theorem~\ref{thrm:EF-relation} follows. $\square$ 

This proof also implies

\begin{corollary} \label{tensor-power}
 Considered as  an 
$\slt$-module with $E$, resp. $F$ acting as $\sum_i[\E_i],$ resp. 
$\sum_i[\F_i],$ the Grothendieck group $K(\O_n)$ 
is isomorphic to the $n$-th  tensor power of the fundamental
two-dimensional representation of $\slt.$
\end{corollary}

\subsubsection {Realization of $\dU$ by projective functors} 
\label{sect-divp} 

Next we provide a realization of the divided powers of 
$E$ and $F$  by projective functors.

Let 
$\E^{(k)}_i$ be the functor from $\O_{i,n-i}$ to $\O_{i+k,n-i-k}$ given 
by tensoring with the $k$-th exterior power of $\Lo$ and then projecting 
onto the submodule with the  generalized central character $\theta_{i+k}:$ 
\begin{equation}
\E^{(k)}_i(M)=\textrm{proj}_{\theta_{i+k}}(\Lambda^k \Lo \o M)
\end{equation}
Similarly, 
\begin{equation}
\F^{(k)}_i:\O_{i,n-i}\to \O_{i-k,n-i+k}
\end{equation} 
is given by 
\begin{equation}
\F^{(k)}_i(M)=\textrm{proj}_{\theta_{i+k}}(\Lambda^k\Ld\o M)
\end{equation}

Denote by $[\E_i^{(k)}], [\F_i^{(k)}]$ the induced homomorphisms of 
the Grothendieck group $K(\O_n)=\oplus_{j\in \Z}K(\O_{j,n-j}).$
Note that $[\E_i^{(k)}],[\F_i^{(k)}]$ map $K(\O_{j,n-j})$ to $0$ unless 
$i=j.$ The following theorem 
is proved in the same way as Theorem~\ref{thrm:EF-relation}.

\begin{theorem} 
\label{thrm:EF-divided-relation} Under the abelian group isomorphism
$$\Upsilon_n: K(\O_n)\to \vn$$ 
the endomorphism $[\E^{(k)}_i],$ resp. $[\F^{(k)}_i]$ 
of $K(\O_n)$ coincides with the endomorphism of $\vn$ given by the action of   
$E^{(k)}1_{2i-n}\in \dU,$ resp. $F^{(k)}1_{2i-n}.$
 In other words, the following 
diagrams are commutative 

\[
\begin{array}{ccc}
K(\O_{n}) & \stackrel{\Upsilon_n}{\longrightarrow} & \vn  \\
\Big\downarrow\vcenter{ \rlap{$[\E^{(k)}_i]$}}     & & 
\Big\downarrow\vcenter{ \rlap{$E^{(k)}1_{2i-n}$}}  \\
K(\O_n) & \stackrel{\Upsilon_n}{\longrightarrow} &  \vn  
\end{array}
\]

\[
\begin{array}{ccc}
K(\O_n) & \stackrel{\Upsilon_n}{\longrightarrow} & \vn  \\
\Big\downarrow\vcenter{ \rlap{$[\F^{(k)}_i]$}}     & & 
\Big\downarrow\vcenter{ \rlap{$F^{(k)}1_{2i-n}$}}  \\
K(\O_n) & \stackrel{\Upsilon_n}{\longrightarrow} &  \vn  
\end{array}
\]

\end{theorem}

Let $S$ be the set $\l E^{(a)}1_j,F^{(a)}1_j | a,j\in \Z\r.$ It 
is a subset of $\dU.$  To each element of $S$ 
we can now associate a projective functor from the category $\O_n$ to 
itself as follows. The functor associated to an element  $x\in S $ will be 
denoted $f_n(x).$ 

\begin{eqnarray*}
& f_n( E^{(a)}1_{j})=& \left\{ 
\begin{array}{ll} 
     \E^{(a)}_{\frac{j+n}{2}} & \mbox{ if $j+n=0(\mbox{mod }2)$} \\
      0                        & \mbox{ otherwise } 
\end{array}
\right.          \\
& f_n( F^{(a)}1_{j})=& \left\{ 
\begin{array}{ll} 
     \F^{(a)}_{\frac{j+n}{2}} & \mbox{ if $j+n=0(\mbox{mod }2)$} \\
      0                        & \mbox{ otherwise } 
\end{array}
\right. 
\end{eqnarray*}

Let $c$ be an arbitrary  product $x_1\dots x_m$ of elements of $S.$ 
To $c$ we associate a functor, denoted 
$f_n(c),$ from $\O_n$ to $\O_n$ by 
$$f_n(c) = f_n(x_1) \circ \dots \circ f_n(x_m).$$

Proposition~\ref{prop:projective-functors-isomorphic} implies 

\begin{theorem}
Let $c_1,\dots, c_s,d_1,\dots, d_t$ be arbitrary products of 
elements of $S.$ 
The endomorphisms of $\vn$ induced by the elements 
 $c_1+\dots +c_s$ and 
$d_1+\dots +d_t$ of $\dU$ coincide  if and only if 
the functors $\oplus_{i=1}^s f_n(c_i)$ and $\oplus_{j=1}^t f_n(d_i)$ 
are isomorphic. 
\end{theorem}

\begin{corollary} There exist isomorphisms of projective functors

\begin{eqnarray*}
\E^{(b)}_{i+a}\circ\E^{(a)}_i  \cong  
(\E^{(a+b)}_i)^{\oplus \binom{a+b}{a}}  & & \\
\F^{(b)}_{i-a}\circ\F^{(a)}_i  \cong 
(\F^{(a+b)}_i)^{\oplus \binom{a+b}{a}} \\ 
\oplus_{k=0}^{min(a,b)}(\E^{(a-k)}_{i-b+k}\circ \F^{(b-k)}_{i}
)^{\oplus\binom{n-i-a+b}{k}} \cong   & & \\ 
\cong \oplus_{l=  0}^{min(a,b)}  
 (\F^{(b-l)}_{i+a-l}\circ \E^{(a-l)}_{i})^{\oplus
\binom{i}{l}} & & 
\end{eqnarray*}
\end{corollary}

\emph{Remark: } We expect that functor isomorphisms in the above 
theorem can be made canonical so that they satisfy certain relations, 
including associativity relations.

\subsubsection{ Canonical basis of $\dU$ and indecomposable 
  projective functors } 
\label{canbas-indecomp}

Recall that the canonical basis $\dB$ of $\dU$ is given by 
\begin{eqnarray*} 
E^{(a)}1_{-i}F^{(b)} & \textrm{ for } & a,b,i\in \N, i\ge a+b \\
F^{(b)}1_i E^{(a)}   & \textrm{ for } & a,b,i\in \N, i > a+b 
\end{eqnarray*}

Conveniently, each element of $\dB$ is a product of 
$E^{(a)}1_i, F^{(a)}1_i$ for various $a$ and $i,$ specifically, 
\begin{eqnarray*} 
E^{(a)}1_{-i}F^{(b)}=E^{(a)}1_{-i}F^{(b)}1_{-i+2b} & 
\textrm{ for } & a,b,i\in \N, i\ge a+b \\
F^{(b)}1_i E^{(a)}=F^{(b)}1_i E^{(a)}1_{i-2a}
   & \textrm{ for } & a,b,i\in \N, i > a+b 
\end{eqnarray*}

In the previous section to each such product and each $n=1,2,\dots $ 
we associated a 
projective functor. Therefore, we can associate a projective functor 
to each element of the canonical basis $\dB$ of $\dU.$ 
To $x\in\dB$ we associate a projective functor $f_n(x)$ 
from $\O_n$ to $\O_n$ by the following rule. 

\begin{eqnarray*} 
f_n(E^{(a)}1_{-i}F^{(b)})& = & \E^{(a)}_{\frac{-i+n}{2}}
\F^{(b)}_{\frac{-i+n}{2}+b}   \\
 & \textrm{ for } & a,b,i\in \N, i\ge a+b, i+n=0(\mbox{mod }2), \\
f_n(F^{(b)}1_i E^{(a)}) &  = & \F^{(b)}_{\frac{i+n}{2}} 
\E^{(a)}_{\frac{i+n}{2}-a} \\
   & \textrm{ for } & a,b,i\in \N, i > a+b, i+n=0(\mbox{mod }2), \\
f_n(E^{(a)}1_{-i}F^{(b)}) & = & f_n(F^{(b)}1_i E^{(a)}) = 0  \\
   & \textrm{ for } & i+n = 1(\mbox{mod }2). 
\end{eqnarray*}

In this way to each canonical basis element $b\in \dB$ there is associated 
an exact functor 
$$f_n(b): \O_n \to \O_n  $$

The multiplication in $\dU$ correspond to composition of projective 
functors: for $x,y\in \dB$ the product $xy$ is 
a linear combination of elements of $\dB$ with 
integral nonnegative coefficients 
 $xy=\sum_{z\in \dU}m_{x,y}^{z}z.$
In turn, the  functor $f_n(x)\circ f_n(y)$ decomposes as a direct sum of 
functors $f_n(z)$ with multiplicities $m_{x,y}^z$: 

$$f_n(x)f_n(y)=\oplus_{z\in \dB}f_n(z)^{\oplus m_{x,y}^z}$$

\begin{theorem}
\label{thrm:indecomposable-proj-functors} Fix $n\in \N.$ 
Let $x\in \dB.$ Then the projective functor $f_n(x)$ is either $0$ or 
isomorphic to an indecomposable projective functor. Moreover, for 
each indecomposable projective functor $A:\O_n\to \O_n$ 
there exists exactly one $x\in \dB$ such that $f_n(x)$ is isomorphic to 
$A.$ 
\end{theorem}

\emph{Proof: }Let $x\in \dB_{2j-n}.$
Recall that in our notations the dominant Verma module in $\O_{j,n-j}$
is $M(1^j0^{n-j}).$  
From the properties of projective functors we know that 
$f_n(x)M(1^j0^{n-j})$ is a projective module in $\O_n,$
$f_n(x)$ is the 
trivial functor if and only if $f_n(x)M(1^j0^{n-j})$ is the trivial 
module, and
$f_n(x)$ is an indecomposable projective functor if and only if 
$f_n(x)M(1^j0^{n-j})$ is an indecomposable projective module. 
Isomorphism classes of projective 
modules in the category of highest weight modules are determined by 
their images in the Grothendieck group. Thus, all computations to 
check whether $f_n(x)M(1^j0^{n-j})$ is indecomposable, trivial, etc. 
can be done in the Grothendieck group of the category $\O_n.$
We claim that 
$[f_n(x)M(1^j0^{n-j})]=0$ or $[f_n(x)M(1^j0^{n-j})]=[P(0^k1^l0^m1^s)]$ 
for some $k,l,m,s$  where, we recall, 
$P(0^k1^l0^m1^s)$ denotes the indecomposable projective cover of 
the simple module with highest weight 
$e_{l+1}+\dots +e_{k+l}+e_{k+l+m+1}+\dots +e_{k+l+m+s}-\rho$ 
(see Section~\ref{section-def-gln}).   

Due to the isomorphism $\Upsilon_n$ between the Grothendieck group 
of $\O_n$ and the abelian group $\vn$ 
and Proposition~\ref{prop:simple-projective-bases} we can
 work with the latest group instead. 
The notations of Section~\ref{sec:canbasis-formulas}
 are used below. 

\begin{prop} 
\begin{eqnarray*}
E^{(a)}1_{i} v(1^b0^c) & = &  \\
 & = &  
\left\{
\begin{array}{ll}  \binom{c}{a} 
(\mbox{Id}^{\otimes b}\o p_c) v(1^{b+a}0^{c-a}) & \mbox{ if }
i=b-c\mbox{ and }c\ge a \\ 
0  & \mbox{otherwise}
\end{array} \right.     \\ 
F^{(a)}1_{i} v(1^b0^c)  & =  &  \\ 
 & =  & \left\{
\begin{array}{ll}  \binom{b}{a} 
(p_b\o \mbox{Id}^{\otimes c}) v(1^{b-a}0^{c+a}) & \mbox{ if }
i=b-c\mbox{ and }b\ge a\\ 
0  & \mbox{otherwise}
\end{array} \right.  
\end{eqnarray*}
\end{prop} 

\emph{Proof:} Clearly, $1_i v(1^b0^c)= v(1^b0^c)$ if $i=b-c$ and 
  $1^iv(1^b0^c)=0$ otherwise. $E^{(a)}$ acts on $v^(1^b0^c)$ in 
  the following way 
 $$ E^{(a)} v(1^b0^c)= \sum_{I} v(1^b I)$$
  where the sum is over all sequences $I$ that contain $a$ ones and 
  $c-a$ zeros. Proposition follows.  $\square$ 

Let $a,b,i$ be non-negative integers with $i\ge a+b.$ 
Then $E^{(a)}1_{-i}F^{(b)}$ is 
an element of the canonical basis $\dB.$ We compute its action on the 
element $v(1^c0^d)$ of $V_1^{\o (c+d)}.$ We restrict to the 
case  $i=2b+d-c$ and $c\ge b$ 
 as otherwise the result is $0.$ The condition $i\ge a+b$ 
can now be written as $c-b\le d-a.$ 

\begin{eqnarray*}
& & E^{(a)}1_{-i}F^{(b)} v(1^c0^d) = \\
& & E^{(a)}\binom{c}{b}(p_c\o Id^{\o d})v(1^{c-b}0^{d+b}) = \\
& & \binom{c}{b}(p_c\o Id^{\o d})E^{(a)}v(1^{c-b}0^{d+b}) = \\
& & \binom{c}{b}(p_c\o Id^{\o d})
 \binom{d+b}{a}(Id^{\o (c-b)}\o p_{d+b})v(1^{c-b+a}0^{d+b-a}) = \\
 & & \binom{c}{b}\binom{d+b}{a}(p_c\o Id^{\o d})
 (Id^{\o (c-b)}\o p_{d+b})v(1^{c-b+a}0^{d+b-a}) = \\
& &  p(0^{b}1^{c-b}0^{d-a}1^{a}) 
\end{eqnarray*}
The last equality follows from formula (\ref{fla:0101}) (substituting 
$j=b, k=c-b, l=d-a, m=a$) and the condition $c-b\le d-a.$ 

In the same fashion, for $a,b,c,d,i\in \Z_+$ such that 
$i=c-d+2a, i\ge a+b$ (which implies $c-b\ge d-a$) 
we get 
\begin{equation}
F^{(b)}1_iE^{(a)} v(1^c0^d) = p(0^b1^{c-b}0^{d-a}1^a)
\end{equation}

Thus, for any $x\in \dB$ and $c,d\in \Z_+$ the element 
$x v(1^c0^d)\in \vn$ is either $0$ or equal to the image under $\Upsilon_n$ 
of the Grothendieck class of the indecomposable projective 
module $P(0^k1^l0^m1^s)$ for some quadruple $(k,l,m,s).$ 

Therefore, for any $x\in \dB_{2j-n}$ and 
a dominant Verma module $M(1^j0^{n-j})$ the projective module 
$f_n(x)M(1^j0^{n-j})$ is 
either the trivial or an indecomposable projective module, i.e. 
the projective functor $f_n(x)$ is either trivial or indecomposable 
projective. All other statements of 
Theorem~\ref{thrm:indecomposable-proj-functors}
 follow easily from the above analysis and the
classification of projective functors (see [BG]). $\square$

\subsubsection{Comultiplication}

\label{sect-comult} 

In previous sections we studied projective functors in categories $\O_n.$ 
We established that on the Grothendieck group level these functors descend 
to the action of generators of $\dU$ on the $n$-th tensor power of 
the fundamental representation $V_1$ and that the composition of functors
descends to the multiplication in $\dU.$ Yet, we do not have a functor
 realization of the whole algebra $\dU$, 
only of its finite-dimensional quotients, also called Schur quotients, that 
are the homomorphic images or $\dU$ in the endomorphisms of $\vn.$  
It is unconvenient to think about comultiplication in $\dU$ if 
only some of its finite quotients are available. We notice, although, that 
comultiplication allows one to introduce a module structure in a  
tensor product and, with a categorification of $\vn$ at hand, we can 
try  to construct functors between  
$\O_n \times \O_m$ and $\O_{n+m}$  corresponding to 
module isomorphism 
\begin{equation}
\label{comultv}
V_1^{\o n}\o V_1^{\o m} \cong V_1^{\o (n+m)}
\end{equation}
Before we considered projective 
projective functors 
$\E_i^{(a)}, \F^{(a)}_i: \O_n \lra \O_n$ for a fixed $n$ and in the notations 
for these functors we suppressed the dependence on $n$. In this section
the rank of $\mf{gl}$ will vary (we'll have functors between 
categories $\O(\gl)\times \O(\mf{gl}_m)$ and $\O(\mf{gl}_{n+m})$ and we 
redenote the functors $\E_i^{(a)}, \F^{(a)}_i:\O_n \lra \O_n $ by 
$\E_{i,n}^{(a)}, \F_{i,n}^{(a)}.$

For the rest of the section we fix positive integers $n$ and $m.$ 
Let $\mf{p}$ be the 
maximal parabolic subalgebra of 
$\mf{gl}_{n+m}$ that contains $\mf{gl}_n\oplus \mf{gl}_m$ and 
$\mf{n}_+,$ i.e. $\mf{p}$ is the standard subalgebra of block uppertriangular 
matrices. Denote by $\mf{u}$ the nilpotent radical of $\mf{p}.$ 
Denote by $\Mp$ the category of finitely-generated 
$U(\mf{p})$-modules that are $\mf{h}$-diagonalizable and 
$U(\mf{n}_+)$-locally nilpotent. 

Let $\Ind$ be the induction functor 
$\Ind:\Mp \to \O(\mf{gl}_{n+m}).$
To a $U(\mf{p})$-module 
$N\in \Mp$ it associates 
the $U(\mf{gl}_{n+m})$-module 
$U(\mf{gl}_{n+m})\otimes_{U(\mf{p})}N.$ 
Recall that we denoted by $F_V$ the functor of tensoring with a finite 
 dimensional $\mf{g}$-module $V.$ 

\begin{lemma} \label{product-and-induction}
Let $V$ be a finite dimensional  $\mf{gl}_{n+m}$-module. Denote by
$F'_V$ the functor from $\Mp$ to $\Mp$ 
given by tensoring with $V,$ considered as a $U(\mf{p})$-module.
Then there is  a canonical isomorphism of functors
\begin{equation} 
F_V\circ \Ind \cong \Ind \circ F'_V
\end{equation} 
\end{lemma}

\emph{Proof: }
For $N\in \Mp$ we have natural isomorphisms
\begin{eqnarray*}
& & \textrm{Hom}_{U(\mf{gl}_{n+m})}(\Ind \circ F'_V(N),F_V\circ \Ind (N)) \\
& = & \textrm{Hom}_{U(\mf{p})}(F'_V(N),F_V\circ \Ind (N)) = \\
& = & \textrm{Hom}_{U(\mf{p})}(V\o N, V\o   (U(\mf{gl}_{n+m})
\otimes_{U(\mf{p})} N )),
\end{eqnarray*}
the first isomorphism coming from the adjointness of induction 
and restriction functors.  
The third  hom-space has a distinguished element coming from the 
$U(\mf{p})$-module map $V \otimes N \lra  V\otimes (U(\mf{gl}_{n+m})
\otimes_{U(\mf{p})} N )$ given  by $v\o n \longmapsto v \o ( 1 \o n).$ 
It is easy to see that the corresponding map of $U(\mf{gl}_{n+m})$ modules 

$$\Ind \circ F'_V(N) \lra F_V \circ \Ind (N)$$ 
is an isomorphism. $\square$ 

Let $Y$ be the one-dimensional representation of $\mf{gl}_n\oplus 
\mf{gl}_m$ which has weight $\frac{m}{2}(e_1+\dots + e_n)$ as 
a representation of $\mf{gl}_n$ and weight $\frac{-n}{2}(e_1+\dots + e_m)$ as 
a representation of $\mf{gl}_m.$ 

Let $\Cfun_0$ be the functor $\O_n\times \O_m\lra \Mp$ defined as 
follows. Tensor a product $M\times N \in \O_n\times \O_m$
with $Y,$ and then make $Y\o M \o N$ 
into a $\mf{p}$-module with the trivial action of the nilpotent radical 
$\mf{u}$ of $\mf{p}.$  

Define $\Cfun: \O_n \times \O_m \lra \O_{n+m}$ to be the composition 
of $\Cfun_0$ and $\Ind:$ 

$$\Cfun = \Ind \circ \Cfun_0.$$

Bifunctor $\Cfun$ is exact and when  
applied to a product of Verma modules 
$M(a_1\dots a_n)\times M(b_1\dots b_m)$ produces the Verma module 
$M(a_1\dots a_n b_1\dots b_m).$ Therefore, we have a commutative diagram 
of isomorphisms of abelian groups:  

$$
\begin{array}{ccc}
  K(\O_n)\times K(\O_m) & \stackrel{\Upsilon_n \times \Upsilon_m}
  {\lra}  &  V_1^{\o n}\times V_1^{\o m} \\
  \downarrow [\Cfun]  &   &    \|    \\
  K(\O_{n+m})  & \stackrel{\Upsilon_{n+m}}{\lra} & V_1^{\o (n+m)} 
\end{array} 
$$

Let us present more evidence that  
 $\Cfun$ is the bifunctor that categorifies the 
identity (\ref{comultv}) by showing how to use $\Cfun$ to 
categorify the comultiplication 
formula $\Delta E= E\o 1+ 1\o E.$

Recall from Section~\ref{functors-EF} that $\Lo$ denotes the $n$-dimensional 
fundamental representation of $\mf{gl}_n.$ We can make $\Lo$ into 
a $U(\mf{p})$-module by making $\mf{u}\oplus \mf{gl}_m\subset \mf{p}$ 
act by $0.$ Similarly, $L_m$ and $L_{n+m}$ are $U(\mf{p})$-modules in 
the natural way and we have a short
 exact sequence of $U(\mf{p})$-modules
\begin{equation}
\label{modules-exact}
0 \to \Lo \to L_{n+m} \to L_m \to 0, 
\end{equation}
which gives rise to a short exact sequence of functors 
from $\Mp$  to $\Mp$: 

$$0\to F_{\Lo}'\to F_{L_{n+m}}'\to F_{L_{m}}'\to 0,$$
where $F_{\Lo}'$ is the functor of tensoring with $\Lo,$ considered 
as a $U(\mf{p})$-module, and so on. 

Functors $\Ind$ and $\Cfun_0$ are exact; and 
composing with them we obtain an  
exact sequence
$$0\to \Ind \circ F_{\Lo}' \circ \Cfun_0 \to \Ind \circ 
F_{L_{n+m}}'\circ \Cfun_0\to 
\Ind \circ F_{L_{m}}'\circ \Cfun_0\to 0$$
of functors. 
We have isomorphisms of functors from $\O_n\times \O_m$ to 
$\Mp:$  
\begin{eqnarray*}
F_{\Lo}' \circ \Cfun_0 & \cong & 
\Cfun_0 \circ (F_{\Lo}\times Id) \\
F_{L_m}' \circ \Cfun_0 & \cong & 
\Cfun_0 \circ (Id \times F_m) 
\end{eqnarray*}
and the isomorphism (see Lemma~\ref{product-and-induction}) 
$$\Ind \circ F_{L_{n+m}}' \cong F_{L_{n+m}} \circ \Ind $$

We thus get an exact sequence of functors 
\begin{equation}
0\to 
\Ind \circ \Cfun_0 \circ (F_{\Lo}\times Id)\to 
F_{L_{n+m}}\circ \Ind \circ \Cfun_0\to 
\Ind \circ \Cfun_0 \circ (Id \times F_{L_m}) 
\to 0,
\end{equation}
and recalling that $\Cfun=\Ind \circ \Cfun_0$ we obtain an exact sequence 
\begin{equation}
\label{exact-seq} 
0\to 
\Cfun \circ (F_{\Lo}\times Id)\to 
F_{L_{n+m}}\circ \Cfun \to 
\Cfun \circ (Id \times F_{L_m}) 
\to 0
\end{equation}
$\E_{i,n}$ is a direct summand of the functor $F_{\Lo}$ 
and from this we derive 

\begin{prop} Exact sequence (\ref{exact-seq}) contains the following exact 
sequence of functors as a direct summand
\begin{equation}
0\to 
\Cfun \circ (\E_{i,n}\times Id)\to 
\E_{i+j,n+m}\circ  \Cfun \to 
\Cfun \circ (Id \times \E_{j,m}) 
\to 0. 
\end{equation}
\end{prop} 
This exact sequence can be considered as a categorification of 
the comultiplication formula 
\begin{equation}
\Delta(E)=E\o 1 +1 \o E
\end{equation}  
In the same fashion, using the exact sequence dual to (\ref{modules-exact})
we obtain an exact sequence of 
functors 
\begin{equation} 
0\to 
\Cfun \circ (Id \times \F_{j,m})\to 
\F_{i+j,n+m}\circ \Cfun \to 
\Cfun \circ (\F_{i,n}\times Id) \to 0
\end{equation} 

We proceed to ``categorify'' the comultiplication rules for 
the divided powers $\E^{(a)}_{i,n}, \F^{(a)}_{i,n}.$ 

The $\mf{gl}_{n+m}$-module $\Lambda^a L_{n+m},$ 
considered as a $\mf{p}$-module, 
admits a filtration 
$\Lambda^a L_{n+m}=G_{a+1}\supset G_{a} \supset \dots \supset G_{0}=0$ 
such that the module $G_{k+1}/G_{k}, k=0,\dots ,a$ is 
isomorphic to $\Lambda^{a-k}\Lo\otimes \Lambda^{k} L_m.$
Therefore, we obtain 

\begin{prop} 
 The exact functor $\E^{(a)}_{i+j,n+m}\circ \Cfun$ has a 
filtration by exact functors 
$$\E^{(a)}_{i+j,n+m}\circ \Cfun=\mc G_{a+1}\supset \mc G_{a} \supset \dots 
\supset \mc G_{0}=0$$
together with short exact sequences of functors
$$0\to \mc G_{k}\to \mc G_{k+1}\to \Cfun \circ 
(\E^{(a-k)}_{i,n}\times \E^{(k)}_{j,m}) \to 0.$$ 

 The exact functor $\F^{(a)}_{i+j,n+m}\circ \Cfun$ has a 
filtration by exact functors 
$$\F^{(a)}_{i+j,n+m}\circ \Cfun=\mc G_{a+1}\supset \mc G_a \supset \dots 
\supset \mc G_{0}=0$$
such that the sequences below are exact
$$0\to \mc G_k\to \mc G_{k+1}\to \Cfun \circ 
(\F^{(k)}_{i,n}\times \F^{(a-k)}_{j,m}) \to 0.$$ 
\end{prop}

This can be considered as a categorification of the comultiplication 
formulas (\ref{comult-div-E}) and (\ref{comult-div-F}). 

\emph{Remark:} 
In [Gr] Grojnowsky uses perverse sheaves to categorify the comultiplication 
rules for $U_q(\mf{sl}_k).$ His approach appears to be Koszul dual to 
ours.

\subsection{Zuckerman functors and Temperley-Lieb algebra} 
\label{TL-and-Z}

\subsubsection{Computations with Zuckerman functors} 
\label{comp-Z}

Let $\mf{g}_i, 1\le i\le n-1$ be a subalgebra of $\mf{gl}_n$ 
consisting of matrices that 
can have non-zero entries only on intersections of $i$-th or $(i+1)$-th 
rows with $i$-th or $(i+1)$-th columns. $\mf{g}_i$ is isomorphic to 
$\mf{gl}_2.$ Denote by $\O_{k,n-k}^i$ the subcategory of $\O_{k,n-k}$ 
consisting of locally $U(\mf{g}_i)$-finite modules. 
Denote by $\Gamma_i: \O_{k,n-k}\lra \O^i_{k,n-k}$ the Zuckerman functor 
of taking the maximal locally $U(\mf{g}_i)$-finite submodule, and by 
$\mc{R}\Gamma_i$ the derived functor of $\Gamma_i.$ 
 The derived functor goes from the bounded derived category 
 $D^b(\O_{k,n-k})$ to $D^b(\O_{k,n-k}^i).$ 
Denote by $\Gamma_i^j$ the cohomology functor 
$\mc{R}^j\Gamma_i: \O_{k,n-k}\to \O_{k,n-k}^j.$ 
This functor is zero if $j>2$ by Proposition~\ref{Zproperties}.

Recall that for a sequence $a_1\dots a_n$ of zeros and ones with exactly 
$k$ ones, the Verma module 
$M(a_1\dots a_n)$ belongs to  $\O_{k,n-k}.$

\begin{prop}
\label{gamma-computed} 
If $a_i=1,a_{i+1}=0,$ 
$$\Gamma^2_i M(a_1\dots a_n)=M(a_1\dots a_n)/
M(a_1\dots a_{i-1}01a_{i+2}\dots a_n)$$
If $a_i=0,a_{i+1}=1,$ 
$$\Gamma^1_i M(a_1\dots a_n)=M(a_1\dots a_{i-1}10a_{i+2}\dots a_n)/
M(a_1\dots a_n)$$
For all other values of $(i,j,a_1,\dots, a_n)$ with  
$i\in {1,\dots , n-1},j\in \Z,a_1,\dots ,a_n\in \l 0,1\r $
we have
$$\Gamma^j_iM(a_1\dots a_n)=0.$$
\end{prop} 

\emph{Proof:} Functors $\Gamma^0_i$ and $\Gamma^2_i$ have  a
simple description as functors of taking the maximal $U(\mf{g}_i)$-locally 
finite submodule/quotient. So the proposition is easy to check for $j\not= 1.$ 
For $j=1$ it is a special case of Proposition 5.5 of [ES]. 
$\square$ 

 The simple module 
$L(a_1\dots a_n)$ with $a_1,\dots, a_n\in \l 0,1\r $ 
and $a_1+\dots +a_n=k$ belongs to $\O^i_{k,n-k}$
if and only if $a_i=1, a_{i+1}=0.$ 
Therefore, a simple object of $\O_{k,n-k}$ cannot belong simultaneously 
to $\O^i_{k,n-k}$ and $\O^{i+1}_{k,n-k}.$ 

\begin{corollary} If a module $M\in \O_{k,n-k}$ lies in both 
$\O^i_{k,n-k}$ and $\O^{i+1}_{k,n-k},$ it is trivial.
\end{corollary}

\begin{prop} For any $M\in \O_{k,n-k}^i, $   
$$\Gamma^j_{i\pm 1}M=0 \mbox{ if }j\not=1$$
\end{prop}

\emph{Proof: } We know that $\Gamma^j_{i\pm 1} 
M$ can be nontrivial only for $0\le j\le 2.$ But $\Gamma^0_{i\pm 1}M$ 
is the maximal locally $\mf{g}_{i\pm 1}$-finite submodule of $M$, 
i.e. the zero module (by the last corollary). 
Similarly, $\Gamma^2_{i\pm 1}M=0.$ 
$\square$

\begin{corollary} 
\begin{enumerate}

\item  Functors $\mc{R}\Gamma_{i\pm 1}[1]$ and 
$\Gamma^1_{i\pm 1}$, restricted to the subcategory 
$\O^i_{k,n-k},$ are naturally equivalent (more precisely, 
the composition of 
$\Gamma^1_{i\pm 1},$ restricted to $\O^i_{k,n-k},$ 
with the embedding functor $\O^{i\pm 1}_{k,n-k}\lra 
D^b(\O_{k,n-k}^{i\pm 1})$ 
is equivalent to the restriction of the functor 
$\mc{R}\Gamma_{i\pm 1}[1]$ to the subcategory $\O_{k,n-k}^i$ of 
$D^b(\O_{k,n-k}^i).$)

\item  The functor 
$\Gamma^1_{i\pm 1}: \O_{k,n-k}^i \to \O_{k,n-k}^{i\pm 1}$
is exact.
\end{enumerate} 
\end{corollary}

Generalized Verma modules for $\mf{g}_i$ are isomorphic to 
the quotients 
$$M(a_1\dots a_{i-1}10a_{i+2}\dots a_n)/
  M(a_1\dots a_{i-1}10a_{i+2}\dots a_n).$$
Denote this quotient module by $M_i(a_1\dots a_{n-2}).$
To simplify notations, if $a_1\dots a_{i-1}a_{i+2}\dots a_n$ is fixed, 
we will denote $M(a_1\dots a_{i-1}10a_{i+2}\dots a_n)$ by $M_{10}$ and 
$M(a_1\dots a_{i-1}01a_{i+2}\dots a_n)$ by $M_{01}.$ Then 
$M_i(a_1\dots a_{n-2})=M_{10}/M_{01}.$

\begin{prop} 
$$\Gamma_{i\pm 1}^1
M_i(a_1\dots a_{n-2})=M_{i\pm 1}(a_1\dots a_{n-2}).$$
\end{prop} 

\emph{Proof: }We have an exact sequence
$$0\to M_{01}\to M_{10}\to M_i\to 0$$
which induces a long exact sequence of Zuckerman cohomology functors
$$\dots \to \Gamma_{i-1}^jM_{01}\to \Gamma_{i-1}^jM_{10}\to 
\Gamma_{i-1}^j M_i\to \Gamma_{i-1}^{j+1}M_{01}\to \dots $$
Consider the following segment of this sequence
\begin{equation}
\label{five-terms}
\Gamma_{i-1}^1M_{01}\to \Gamma_{i-1}^1M_{10}
\to \Gamma_{i-1}^1M_i\to \Gamma_{i-1}^2 M_{01}
\to \Gamma_{i-1}^2 M_{10} 
\end{equation}
From Proposition~\ref{gamma-computed}, we find that 
$$\Gamma_{i-1}^1M_{01}=\Gamma_{i-1}^2M_{10}=0$$
and (\ref{five-terms}) becomes a short exact sequence
$$0 \to \Gamma_{i-1}^1M_{10}
\to \Gamma_{i-1}^1M_i\to \Gamma_{i-1}^2 M_{01}
\to 0 $$
We proceed by considering two different cases: 

(i) If $a_{i-1}=0,$ then by Proposition~\ref{gamma-computed}
we have $\Gamma_{i-1}^1M_{10}=M_{i-1}(a_1\dots a_{n-2}) $ and 
$\Gamma_{i-1}^2 M_{01}=0.$ 

(ii) If $a_{i-1}=1,$ then 
$\Gamma_{i-1}^2M_{01}=M_{i-1}(a_1\dots a_{n-2}) $ and 
$\Gamma_{i-1}^1 M_{10}=0.$ 

In both cases we get $\Gamma_{i-1}^1M_i=
M_{i-1}(a_1\dots a_{n-2}).$ The proof for $\Gamma_{i+1}$ is 
the same. 
$\square$

\subsubsection{Equivalences of categories} 
\label{equiv-of-cat}

Let $\varepsilon_i$ be the inclusion functor 
$\O^i_{k,n-k}\to \O_{k,n-k}.$
Consider a pair of functors 

\begin{eqnarray*}
\Gamma_{i-1}^1\varepsilon_i & : 
\O^i_{k,n-k}  & \to \O^{i-1}_{k,n-k} \\
\Gamma_{i}^1\varepsilon_{i-1} & : 
\O^{i-1}_{k,n-k}  & \to \O^{i}_{k,n-k} 
\end{eqnarray*}

These two functors are exact, take generalized Verma modules to 
generalized Verma modules and the compositions 
$\Gamma_{i-1}^1 \varepsilon_i
\Gamma_{i}^1 \varepsilon_{i-1}$ and 
$\Gamma_{i}^1\varepsilon_{i-1}
\Gamma_{i-1}^1\varepsilon_i$ are identities on generalized Verma 
modules.  

\begin{theorem} 
\label{functors-and-isotopies}
Functors $\Gamma_{i - 1}^1\varepsilon_i$ 
and $\Gamma_i^1\varepsilon_{i-1}$ 
are equivalences of 
categories $\O_{k,n-k}^i$ and $\O_{k,n-k}^{i- 1}.$
The composition $\Gamma_{i}^1\varepsilon_{i- 1}
\Gamma_{i- 1}^1\varepsilon_{i}$ is  isomorphic
 to the identity functor from the category $\O_{k,n-k}^i$ to itself. 
\end{theorem} 

\emph{Proof:} By Lemma~\ref{category-equivalence} 
it remains to show that functors 
$\Gamma_{i-1}^1\varepsilon_i$ and 
$\Gamma_{i}^1\varepsilon_{i-1}$ are two-sided adjoint. 
We recall that $\Gamma_i:\O_{k,n-k}\to \O^i_{k,n-k}$ is 
right adjoint to the inclusion functor $\varepsilon_i.$ 
Besides, in the derived category, 
the derived functor $R\Gamma_i$ is isomorphic to the left adjoint of 
the shifted inclusion functor $\varepsilon_i[2].$

Let $M,N\in \O^i_{k,n-k}.$ 
We have natural vector space isomorphisms 
\begin{eqnarray*}
\mbox{Hom}((\Gamma_{i-1}^1\varepsilon_i )M , N ) & \cong & 
\mbox{Hom}((R\Gamma_{i-1}[1]\circ \varepsilon_i) M , N ) 
\\
& \cong & 
\mbox{Hom}(\varepsilon_i M , \varepsilon_{i-1}N[1] ) 
\\
& \cong & 
\mbox{Hom}(M , (R\Gamma_{i})\varepsilon_{i-1}N[1] ) 
\\
& \cong & 
\mbox{Hom}(M , (R\Gamma_{i})[1]\varepsilon_{i-1}N) 
\\
& \cong & 
\mbox{Hom}(M , \Gamma_{i}^1\varepsilon_{i-1} N) 
\end{eqnarray*}
which imply that $\Gamma_{i-1}^1\varepsilon_i$  
 is left adjoint to $\Gamma_i^1\varepsilon_{i-1}.$ 
A similar computation tells us that  $\Gamma_{i-1}^1\varepsilon_i$  
 is also right  adjoint to $\Gamma_i^1\varepsilon_{i-1}.$ 

\emph{Remark:} For an abelian category $\mc{A}$ the embedding functor 
$\mc{A}\lra D^b(\mc{A})$ is fully faithful and for $M,N\in \mc{A}$ 
we have canonical isomorphisms of hom-spaces 
$$\mbox{Hom}_{\mc{A}}(M,N) = \mbox{Hom}_{D^b(\mc{A})}(M,N).$$ 
That permits us in the chain of identities above to go freely between 
hom spaces in $\O^i_{k,n-k}$ and $D^b(\O^i_{k,n-k}).$ $\square$

\begin{corollary} The categories $\O_{k,n-k}^i$ and $\O_{k,n-k}^j$ 
are equivalent for all $i,j, 1\le i,j \le n-1.$
\end{corollary}

\begin{prop} 
\label{opiat-dvoika}
Let $M\in \O_{k,n-k}^i.$ Then 
$$\Gamma_i^j\varepsilon_iM=
\left\{ 
\begin{array}{ll} 
M & \mbox{ if }j=0\mbox{ or }j=2 \\
0 & \mbox{ otherwise }
\end{array} \right. 
$$
\end{prop}

\emph{Proof: } $\Gamma_i^0\varepsilon_iM$ 
is the maximal $\mf{g}_i$-locally finite 
submodule of $M$ and $\Gamma_i^2\varepsilon_iM$ 
is the maximal $\mf{g}_i$-locally 
finite quotient of $M.$ Thus, $\Gamma_i^0\varepsilon_iM=
\Gamma_i^2\varepsilon_iM=M.$ 
It remains to check that $\Gamma_i^1\varepsilon_iM=0$ 
for any $M\in \O_{k,n-k}^i.$ 

The derived functor $\mc{R}\Gamma_i$ is exact and induces a map of 
Grothendieck groups 
$$[\mc{R}\Gamma_i]: K(\O_{k,n-k})\to K(\O_{k,n-k}^i).$$
Using Proposition~\ref{gamma-computed} 
we can easily compute $[\mc{R}\Gamma_i].$ 
The bases of $K(\O_{k,n-k})$ and 
$K(\O_{k,n-k}^i)$ consisting of images of Verma modules, respectively
generalized Verma modules, are convenient for writing down $[R\Gamma_i]$ 
explicitly: 

\begin{eqnarray*} 
& [\mc{R}\Gamma_i(M(a_1\dots a_{i-1}00 a_i\dots a_{n-2}))]= & 
0\\
& [\mc{R}\Gamma_i(M(a_1\dots a_{i-1}11 a_i\dots a_{n-2}))]= & 
0 \\
& [\mc{R}\Gamma_i(M(a_1\dots a_{i-1}10 a_i\dots a_{n-2}))]= & 
[M_i(a_1\dots a_{n-2})] \\
& [\mc{R}\Gamma_i(M(a_1\dots a_{i-1}01 a_i\dots a_{n-2}))]= & 
-[M_i(a_1\dots a_{n-2})] 
\end{eqnarray*}
where, we recall
$$M_i(a_1\dots a_{n-2})=M(a_1\dots a_{i-1}10a_i\dots a_{n-2})/
M(a_1\dots a_{i-1}01a_i\dots a_{n-2}).$$

Therefore, 
$$[\mc{R}\Gamma_i\circ \varepsilon_i (M_i(a_1\dots a_{n-2}))]=
[M_i(a_1\dots a_{n-2})]\oplus
[M_i(a_1\dots a_{n-2})]$$
and 
$$[(\mc{R}\Gamma_i\circ\varepsilon_i) M]=[M]\oplus [M]$$
for any $M\in \O^i_{k,n-k}.$ On the other hand, 
$$[(\mc{R}\Gamma_i\circ\varepsilon_i)M]=
[\Gamma_i^0\varepsilon_iM]-[\Gamma_i^1\varepsilon_iM]+
[\Gamma_i^2\varepsilon_iM]=[M]-[\Gamma_i^1\varepsilon_iM]+[M].$$
Thus, $[\Gamma_i^1\varepsilon_iM]=0$ for any $M\in \O^i_{k,n-k}$
 and, hence, $\Gamma_i^1\varepsilon_iM=0$ for any $M\in \O^i_{k,n-k}.$ 
$\square$

\begin{prop} 
\label{id-and-id}
Restricting to the subcategory 
$D^b(\O_{k,n-k}^i),$ 
we have an equivalence of functors 
\begin{equation}
\label{double-id}
\mc{R}\Gamma_i\circ\varepsilon_i\cong Id\oplus Id[-2].
\end{equation}
\end{prop} 

\emph{Proof: }Consider the natural transformation 
$$e:Id \to \mc{R}\Gamma_i\circ \varepsilon_i.$$
coming from the adjointness of $\mc{R}\Gamma_i$ and $\varepsilon_i.$ 
Let $C_i$ be the functor which is the cone of $e.$ Then 
by [MP], Lemma 3.2,  we have an isomorphism 
of functors 
$$\mc{R}\Gamma_i\circ \varepsilon_i=Id\oplus C_i.$$
From Proposition~\ref{opiat-dvoika}
$$C_i^jM=\left\{
\begin{array}{ll} M  & \mbox{ if }j=2 \\
                  0  & \mbox{ if }j\not= 2 
\end{array}  \right. 
$$
Therefore, $C_i=Id[-2]$ and we have the isomorphism (\ref{double-id}).

$\square$ 

\subsubsection{A realization of the Temperley-Lieb algebra by functors} 
\label{TL-by-functors} 

Define functors $\mc{V}_i, 1\le i\le n-1$ from 
$D^b(\O_{k,n-k})$ to $D^b(\O_{k,n-k})$ by 
\begin{equation}
\mc{V}_i = \varepsilon_i \circ \mc{R}\Gamma_i [1]
\end{equation} 

\begin{theorem}\label{TL-via-derived} 
There are natural equivalences of functors 
\begin{eqnarray}
\label{dva-i-dva}
  (\mc{V}_i )^2 & \cong & \mc{V}_i [-1] \oplus \mc{V}_i [1] \\
\label{commutativity-from-afar}
   \mc{V}_i \mc{V}_j & \cong & \mc{V}_j\mc{V}_i \mbox{ for }|i-j|>1 \\
\label{little-braiding}
   \mc{V}_i \mc{V}_{i\pm 1}\mc{V}_i & \cong & \mc{V}_i 
\end{eqnarray}
\end{theorem}  

\emph{Proof:} Isomorphism (\ref{dva-i-dva}) follows from 
Proposition~\ref{id-and-id}. Isomorphism (\ref{commutativity-from-afar})
is implied by a commutativity isomorphism $\Gamma_i\Gamma_j \cong 
\Gamma_j \Gamma_i$ for $|i-j|>1.$ 
The last isomorphism is a corollary of Theorem~\ref{functors-and-isotopies}. 
$\square$ 

Summing over all $k$ from $0$ to $n$ we obtain functors 
$$\mc{V}_i: D^b(\O_n) \lra D^b(\O_n)$$ 
together with isomorphisms 
(\ref{dva-i-dva})-(\ref{little-braiding}). On the Grothendieck 
group level these functors descend to the action of the Temperley-Lieb 
algebra $TL_{n,1}$ on $\vn.$ We thus have a categorification of 
the action of the Temperley-Lieb algebra on the tensor product $\vn.$ 
This action is faithful, so we can loosely say that we have 
a categorification of the Temperley-Lieb algebra $TL_{n,1}$ itself.  
In fact if we consider the shift by $1$ in the derived category as 
the analogue of the multiplication by $q$, then we have categorified 
the Temperley-Lieb algebra $TL_{n,q}.$

Recall from Section~\ref{TL-algebra} that the element $U_i$ of the 
Temperley-Lieb algebra is the product of morphisms $\cup_{i,n-2}$ and 
 $\cap_{i,n}$ of the Temperley-Lieb category.  Morphisms 
 $\cap_{i,n}$ and $\cup_{i,n-2}$  go  between objects $\overline{n}$ 
 and $\overline{n-2}$ of the Temperley-Lieb category.  
  On the other hand, 
 the functor $\mc{V}_i,$ which categorifies $U_i,$ is the composition 
  of functors $\varepsilon_i$  and $ \mc{R}\Gamma_i [1]$. 
 We now explain how to modify the latter functors into functors 
  between derived categories $D^b(\O_n)$ and $D^b(\O_{n-2})$ that 
  cat be viewed as categorifications of morphisms $\cap_{i,n}$ and 
  $\cup_{i,n-2}.$ 

Let $\varsigma_{n} $ be the functor from 
$\O_{k-1,n-k-1}$ to $\O_{k,n-k}^1$ given as follows. 
First we tensor an $M\in \O_{k-1,n-k-1}$ with the fundamental 
representation $L_2$ of $\mf{gl}_2$ to get a 
$\mf{gl}_2\oplus \mf{gl}_{n-2}$-module $L_2\otimes M.$ 
Let $Y$ be the one-dimensional $\mf{gl}_2\oplus \mf{gl}_{n-2}$-module 
with weight $\frac{n-3}{2}(e_1+e_2)$ relative to $\mf{gl}_2$ and 
$-(e_1+\dots +e_{n-2})$ relative to $\mf{gl}_{n-2}.$ 
Let $\mf{p}$ be the maximal parabolic subalgebra of $\mf{gl}_n$ that 
contains $\mf{gl}_2\oplus \mf{gl}_{n-2}$ and the subalgebra of 
upper triangular matrices. Then $Y\o (L_2\otimes M)$ is naturally a 
$\mf{p}$-module with the nilradical of $\mf{p}$ 
acting trivially. 
Now we parabolically induce from $\mf{p}$ to $\mf{gl}_n$ and define
$$\varsigma_n(M)=U(\mf{gl}_n)\otimes_{U(\mf{p})}(Y\o L_2\o M).$$

Let $\nu_n$ be the functor 
from $\O_{k,n-k}^1$ to $\O_{n-1,n-k-1}$ defined as follows. 
For an $M\in \O_{k,n-k}^1,$ take the sum of the weight subspaces
of $M$ of weights $e_1+x_3e_3+\dots +x_ne_n-\rho_n,$ where 
$x_3,\dots ,x_n\in \Z$ and $\rho_n$ is the half-sum of the positive roots 
of $\mf{gl}_n.$ This direct sum is a $\mf{gl}_{n-2}$-module in a natural way. 
Define $\nu_n(M)$ as the tensor product of this module with 
the one-dimensional $\mf{gl}_{n-2}$-module of weight $e_1+\dots +e_{n-2}.$  

\begin{prop} Functors $\varsigma_n$ and $\nu_n$ are mutually inverse 
equivalences of categories $\O_{k-1,n-k-1}$ and $\O_{k,n-k}^1.$ 
\end{prop} 

We omit the proof as it is quite standard. 

\begin{corollary} The categories $\O_{k,n-k}^i,1\le i\le n-1$ and 
$\O_{k-1,n-k-1}$ are equivalent. 
\end{corollary} 

Denote by $\Xi_{n,i}$ the equivalence of categories 
$$\Xi_{n,i}: \O_{k,n-k}^i \lra \O_{k-1,n-k-1}$$
given by the composition 
$$\Xi_{n,i}= \nu_n \circ \Gamma^1_1 \circ \varepsilon_2 
\circ \Gamma^1_2 \circ \dots \varepsilon_{i-1}
\circ \Gamma_{i-1}^1\circ \varepsilon_i,$$
 i.e., $\Xi_{n,i}$ is the composition of equivalences
of categories 
$$ \O_{k,n-k}^i \stackrel{\cong}{\lra} \O_{k,n-k}^{i-1} 
\stackrel{\cong}{\lra} \dots \stackrel{\cong}{\lra} 
 \O_{k,n-k}^1 \stackrel{\cong}{\lra} \O_{k-1,n-k-1}. $$

Denote by $\Pi_{n,i}$ the equivalence of categories 
$$\Pi_{n,i}: \O_{k,n-k} \lra \O_{k+1,n+1-k}^i$$
given by 
$$\Pi_{n,i}= \Gamma^1_i\circ \varepsilon_{i-1} 
\circ \Gamma^1_{i-1} \circ \dots \circ \varepsilon_2\circ \Gamma^1_2
\circ \varepsilon_1\circ\varsigma_n .$$

Denote the derived functors
of these functors  by  $\mc{R}\Xi_{n,i}$ and $\mc{R}\Pi_{n,i}$.
Define functors 
\begin{eqnarray*} 
\cap_{i,n}: &  D^b(\O_{k,n-k}) & \lra D^b(\O_{k-1,n-k-1}) \\
\cup_{i,n}: &  D^b(\O_{k,n-k}) & \lra D^b(\O_{k+1,n+1-k})  
\end{eqnarray*} 
by 
\begin{eqnarray}
\label{here-is-cap}
\cap_{i,n} & = & \mc{R}\Xi_{n,i}\circ \mc{R}\Gamma_i[1] \\
\label{here-is-cup}
\cup_{i,n} & = & \varepsilon_i \circ \mc{R}\Pi_{n,i}  
\end{eqnarray} 
Recall from Section~\ref{TL-algebra} defining relations
(\ref{first-relation})-(\ref{last-relation}) for the Temperley-Lieb 
category $TL.$ 

\begin{conjecture} There are natural equivalences 
(\ref{first-relation})-(\ref{sixth-relation})
 with functors $\cap_{i,n}$ and $\cup_{i,n}$ defined by 
(\ref{here-is-cap}) and (\ref{here-is-cup}). 
\end{conjecture}

The first two of these equivalences follow from the results of this section. 
The relation (\ref{last-relation}) will become 
$\cap_{i,n+2}\circ \cup_{i,n}\cong Id[1]\oplus Id[-1].$

We next state a conjecture on a functor realization of the category of 
 tangles in $\R^3.$ Consider the following 3 elementary tangles 

\begin{center} \epsfig{file=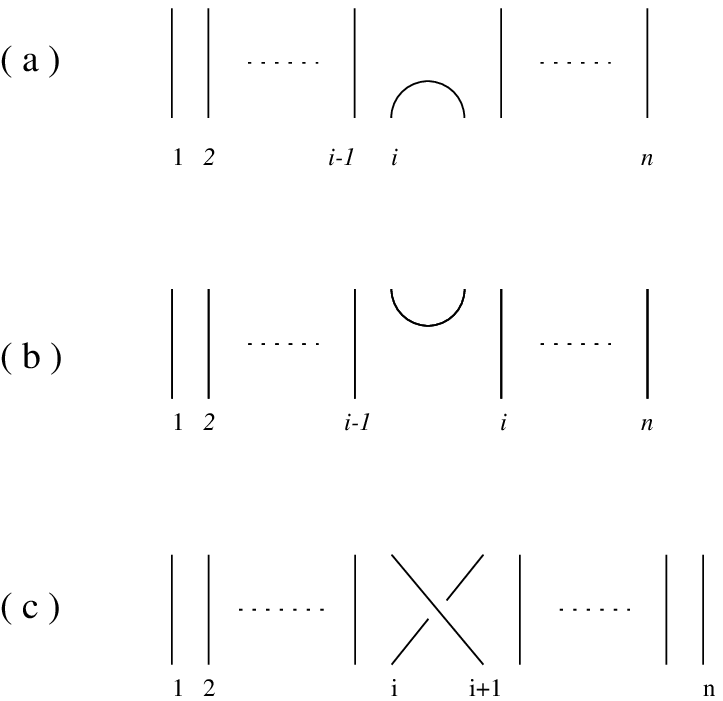} \end{center}

Every tangle in 3-space can be presented as a concatenation of these 
elementary tangles. We associate to these 3 types of tangles the 
following functors: 

To the 
tangle (a) associate functor $\cap_{i,n}$ given by the formula (\ref{here-is-cap}). 

To the tangle (b) associate functor $\cup_{i,n}$ given by the formula 
(\ref{here-is-cup}). 

To the tangle (c) associate functor $R_{i,n}$ from $D^b(\O_n)$ to $D^b(\O_n)$ 
which is the cone of the adjointness morphism of functors 
$$ \varepsilon_i \circ \mc{R}\Gamma_i \lra Id. $$
Given a presentation $\alpha$ of a tangle $t$ as a composition 
of elementary tangles of types (a)-(c), to $\alpha$ we associate the 
functor $f(\alpha)$ which is the corresponding composition of 
functors $\cap_{i,n}, \cup_{i,n}, R_{i,n}.$ 

\begin{conjecture}
\label{tangle-invariants-sing}
 Given two such presentations $\alpha, \beta$ of 
a tangle $t$, functors $f(\alpha)$ and $f(\beta)$ are 
isomorphic, up to shifts in the derived category. 
\end{conjecture}

This conjecture, if true, will give us functor invariants of tangles. 
Certain natural  transformations between these 
functors, corresponding to adjointness morphisms between 
the (shifted) identity functor and compositions 
of $\varepsilon_i $ and $\mc{R}\Gamma_i,$ are expected to produce 
invariants of 2-tangles and 2-knots (to be discussed elsewhere). 
When we have a link rather than a tangle, the associated functor 
goes between categories of complexes of vector spaces up to homotopies, 
and the resulting invariants of links will be $\Z$-graded homology groups.


\section{Parabolic categories} 
\label{parab-section}
\subsection{Temperley-Lieb algebra and projective functors}
\label{parab-tl} 
\subsubsection{On and off the wall translation functors in parabolic 
 categories} 
\label{on-and-off} 

Let $\mu$ be an integral dominant regular weight and $\mu_i, i=1,\dots , 
n-1$ an integral dominant subregular weight on the $i$-th wall. 
Let $\O_{\mu}$ and $ \O_{\mu_i}$
be the subcategories of $\O(\gl)$ of modules with generalized 
central characters $\eta(\mu)$ and $\eta(\mu_i).$ 
Then $\O_{\mu}$ is a regular block of $\O(\gl)$ and $\O_{\mu_i}$ is a 
subregular block of $\O(\gl).$ 
Verma modules $M_{\mu}$ and $M_{\mu_i}$ with highest weights 
$\mu-\rho$ and $\mu_i-\rho$ are dominant Verma 
modules in the corresponding categories. 

Let $\Tui, \Ti$ be translation functors 
on and off the $i$-th wall

\begin{eqnarray*}
\Tui : & \O_{\mu} \lra & \O_{\mu_i} \\
\Ti : & \O_{\mu_i} \lra & \O_{\mu} 
\end{eqnarray*}
These functors are defined up to an isomorphism by the condition that 
they are projective functors between $\O_{\mu}$ and 
$\O_{\mu_i}$ and 
\begin{enumerate}
\item Functor $\Tui$ takes the Verma module $M_{\mu}$ to the Verma module 
$M_{\mu_i}.$ 
\item Functor $\Ti$ takes the Verma module $M_{\mu_i}$ to the projective 
module $P_{s_i\mu}$ where $s_i$ is 
the transposition $(i,i+1).$ On the Grothendieck group level, 
$$[\Ti M_{\mu_i}] = [M_{\mu}] +  [M_{s_i\mu}].$$
\end{enumerate} 

Let $\mf{p}_k$ be the maximal parabolic subalgebra of $\gl$ 
such that $\mf{p}_k\supset \mf{n}_+\oplus \mf{h}$ and 
the reductive subalgebra of $\mf{p}_k$ is $\mf{gl}_k\oplus \mf{gl}_{n-k}.$
Let 
$\O^{k,n-k}$, resp $\O^{k,n-k}_i$ 
be the full subcategory of $\O_{\mu}$, resp. $\O_{\mu_i}$ 
consisting 
of modules that are locally $U(\mathfrak{p}_k)$-finite. 

From now on we fix $k$ between $0$ and $n$. Let $\tau_i^{i + 1}$ 
be the composition of $\Ti$ and $T^{i+1}$: 
$$\tau_i^{i+1}= T^{i+1}\circ \Ti.$$
This is a functor from $\O_{\mu_i}$ to $\O_{\mu_{i+1}}.$ 
Similarly, let $\tau_i^{i-1}$ be the functor from 
$\O_{\mu_i}$ to $\O_{\mu_{i-1}}$ given by 
$$\tau_{i}^{i-1}= T^{i-1} \circ \Ti.$$
 Projective functors preserve subcategories of $U(\mf{p}_k)$-locally 
finite modules and thus functors $\tau_i^{i \pm 1}$ restrict to 
functors from $\O_i^{k,n-k}$ to $ \O_{i \pm 1}^{k,n-k}.$ 
Our categorification of the Temperley-Lieb algebra by projective functors 
is based on the following beautiful result of Enright and Shelton:  

\begin{theorem}\label{project-equivalence} 
Functors $\tau_i^{i \pm 1}$ establish equivalences 
of categories $\O_i^{k,n-k}$ and $\O_{i\pm 1}^{k,n-k}.$ 
\end{theorem}

\emph{Proof:} See [ES], Lemma 10.1 for a proof of a 
slightly more general statement. 
An alternative proof that we give below uses  

\begin{lemma} The functor $\tau_{i+1}^i \tau_i^{i+1}$ restricted to 
$\O_i^{k,n-k}$ is isomorphic to the identity functor.  
\end{lemma} 

\emph{Proof} We first study this functor as a projective functor from 
the subregular block $\O_{\mu_i}$ to itself. We will  show that this 
functor is a direct sum of the identity functor and another projective 
functor that vanishes when restricted to $\O_i^{k,n-k}.$ 
 
An isomorphism class of a projective 
functor is determined by its action on the dominant Verma module. 
So let us compute the action of $\tau_{i+1}^i\tau_i^{i+1}$ 
on $M_{\mu_i}$ on the Grothendieck group level. 

\begin{eqnarray*}
 [\tau_{i+1}^i \tau_i^{i+1} M_{\mu_i}] & = & 
  [\Tui T_{i+1} T^{i+1} \Ti M_{\mu_i}] \\ 
                                     & = & 
  [\Tui T_{i+1} T^{i+1} (M_{\mu}\oplus M_{s_i \mu}) ] \\
                                     & = & 
  [\Tui T_{i+1} (M_{\mu_{i+1}} \oplus M_{s_i\mu_{i+1}}) ] \\
                                     & = & 
  [\Tui ( M_{\mu} \oplus M_{s_{i+1}\mu} \oplus M_{s_i \mu} 
   \oplus M_{s_is_{i+1}\mu}]  \\
                                     & = & 
  [M_{\mu_i} \oplus M_{s_{i+1}\mu_i} \oplus M_{s_i \mu_i} 
   \oplus M_{s_is_{i+1}\mu_i}]  \\
                                     & = & 
  [M_{\mu_i}] +  [M_{s_{i+1}\mu_i}] +  [M_{s_i \mu_i}] + 
   [M_{s_is_{i+1}\mu_i}]  \\
                                     & = & 
  [M_{\mu_i}] +  [M_{s_{i+1}\mu_i}] +  [M_{\mu_i}] + 
   [M_{s_is_{i+1}\mu_i}] 
\end{eqnarray*} 
                                      
The projective module $P_{s_is_{i+1}\mu_i}$ 
decomposes in the Grothendieck group as the following sum 
$$[P_{s_is_{i+1}\mu_i}]=  
[M_{\mu_i}] + [M_{s_{i+1}\mu_i}] + [M_{s_is_{i+1}\mu_i}].$$
Therefore,  
$$[\tau_{i+1}^i \tau_i^{i+1} M_{\mu_i}] = [M_{\mu_i}] + 
[P_{s_is_{i+1}\mu_i}].$$
From the classification of projective functors (see [BG]), we 
derive that $\tau_{i+1}^i \tau_i^{i+1}$ is isomorphic to 
the direct sum of the identity functor and an indecomposable 
 projective functor that 
takes $M_{\mu_i}$ to the indecomposable projective module 
$P_{s_is_{i+1}\mu_i}.$ 
Denote this functor by $\wp.$ Then 
$$\tau_{i+1}^i \tau_i^{i+1}\cong Id \oplus \wp.$$

Let $\xi_i$ be an integral dominant weight on the intersection of 
the $i$-th and $(i+1)$-th walls. We require that $\xi_i$ be 
a generic weight with these conditions, i.e. $\xi_i$ does not 
lie on any other walls. 
Let $\O_{\xi_i}$ be the subcategory of $\O(\gl)$ consisting 
of modules with generalized central character $\eta(\xi_i).$ 

Let $T_{\mu_i}^{\xi_i}$ and $T_{\xi_i}^{\mu_i}$ be translation 
functors from $\O_{\mu_i}$ to $\O_{\xi_i}$ and back. Then 
$T_{\mu_i}^{\xi_i}$ takes the Verma module $M_{\mu_i}$ to the 
Verma module $M_{\xi_i}$ while $T_{\xi_i}^{\mu_i}$ takes 
$M_{\xi_i}$ to $P_{s_is_{i+1}\mu_i}.$ 
Therefore, functor $\wp $ is isomorphic to the composition 
$T_{\xi_i}^{\mu_i}T_{\mu_i}^{\xi_i}$ and 
$$\tau_{i+1}^i \tau_i^{i+1}\cong Id \oplus 
T_{\xi_i}^{\mu_i}T_{\mu_i}^{\xi_i}.$$ 
 The category $\O_{\xi_i}$ contains no $U(\mf{p}_k)$-locally finite 
modules other than the zero module. Hence, the functor $T_{\mu_i}^{\xi_i},$ 
restricted to the subcategory $\O_i^{k,n-k}$ is the zero functor. 
Therefore, $\tau_{i+1}^i\tau_i^{i+1}$ , restricted to $\O_i^{k,n-k},$ 
is isomorphic to the identity functor. This proves the lemma. 
$\square$

In exactly the same fashion we establish that $\tau_i^{i+1}\tau_{i+1}^i$ 
is isomorphic to the identity functor from $\O_{i+1}^{k,n-k}$ to itself. 
Therefore functors  $\tau_i^{i+1}$ and $\tau_{i+1}^i$ are mutually 
inverse and provide an equivalence of categories $\O_{i}^{k,n-k}$ and 
$\O_{i+1}^{k,n-k}.$ 
$\square$ 

\subsubsection{Projective functor realization of the Temperley-Lieb  
  algebra} 
\label{pr-and-TL} 

Define the  functor $\CU_i, i=1,\dots n-1$ from $\O^{k,n-k}$ to $\O^{k,n-k}$ 
as the composition of functors $\Ti$ and $\Tui$: 
$$\CU_i=\Ti \circ \Tui$$

\begin{prop}
\label{TL-projective}
There are equivalences of functors
\begin{eqnarray*}
\CU_i\circ \CU_j   & \cong & \CU_j\circ \CU_i\mbox{ for }|i-j|>1 \\
\CU_i\circ \CU_i   & \cong & \CU_i\oplus \CU_i \\
\CU_i\circ \CU_{i\pm 1}\circ \CU_i & \cong &  \CU_i 
\end{eqnarray*}
\end{prop}

\emph{Proof: }The first two equivalences hold even if we consider 
$\CU_i$ as the functor in the bigger category $\O_{\mu}.$ For example, 
$$\CU_i\circ \CU_i= \Ti\circ \Tui\circ \Ti \circ \Tui= 
(\Ti \circ \Tui) \oplus (\Ti\circ \Tui) =
\CU_i \oplus \CU_i$$
where the second equality follows from the result that the 
composition $\Tui\circ \Ti$ of projective functors off and on the wall 
is the direct sum of two copies of the identity functor. 

For the third equivalence the restriction to $\O^{k,n-k}$ is 
absolutely necessary. Then 
\begin{eqnarray*} 
\CU_i \CU_{i+1}\CU_i & = & \Ti\Tui T_{i+1}T^{i+1}\Ti \Tui \\
                     & = & \Ti \tau_{i+1}^i \tau_i^{i+1}\Tui  \\
                     & = & \Ti \Tui  \\
                     & = & \CU_i 
\end{eqnarray*} 
$\square$

  This proposition gives a functor realization of the Temperley-Lieb 
  algebra by projective functors in parabolic categories $\O^{k,n-k}.$ 
Suitable products of generators $U_i$ produce a basis of the Temperley-Lieb 
 algebra $TL_{n,q},$ which admits a graphical interpretation: elements 
  of this basis correspond to isotopy classes of systems of $n$ simple 
 disjoint arcs in the plane connecting $n$ bottom and $n$ top points. 
 Moreover, this basis is related (see [FG]) 
 to the Kazhdan-Lusztig basis in the Hecke algebra as well as to Lusztig's 
  bases in tensor products of $U_q(\mf{sl}_2)$-representations (see [FK]). 
 Proposition~\ref{TL-projective} implies that to an element of this basis, 
  we can associate a projective functor from $\O^{k,n-k}$ to $\O^{k,n-k},$
  which is defined as a suitable composition of functors 
 $\CU_i, 1\le i\le n-1.$  
  We conjecture that these compositions of $\CU_i$'s are 
  indecomposable, and, in turn, an indecomposable projective functor 
  from $\O^{k,n-k}$ to $\O^{k,n-k}$ 
  is isomorphic to one of these compositions
 (compare with Theorem~\ref{thrm:indecomposable-proj-functors}).

  In [ES] Enright and Shelton, among other things, constructed an 
equivalence of categories $\O_1^{k,n-k}$ and $\O^{k-1,n-k-1}$ 
(see [ES], $\S 11$). This equivalence allows us to factorize 
$\CU_i$ as a composition of a functor from $\O_{k,n-k}$ to 
$\O_{k-1,n-k-1}$ and a functor from $\O_{k-1,n-k-1}$ to $\O_{k,n-k}$. 
This is very much in line with the factorization of the element $U_i$ 
of the Temperley-Lieb algebra as the composition  
$\cup_{i,n-2}\circ \cap_{i,n}$ 
of morphisms $\cup_{i,n-2}$ and $\cap_{i,n}$ of the Temperley-Lieb category. 

We now offer the reader a conjecture on 
realizing the Temperley-Lieb category via functors between parabolic 
categories $\O^{k,n-k}.$ Let 
\begin{equation}
\zeta_n:\O^{k,n-k}_1 \stackrel{\cong}{\lra} \O^{k-1,n-k-1}
\end{equation}
be the Enright-Shelton equivalence of categories. 
Introduce functors 
\begin{eqnarray*} 
\cap_{i,n}: &  \O^{k,n-k} & \lra \O^{k-1,n-k-1} \\
\cup_{i,n}: &  \O^{k,n-k} & \lra \O^{k+1,n+1-k} 
\end{eqnarray*} 
given by 
\begin{eqnarray}
\label{another-cap}
\cap_{i,n} & = & \zeta_n \circ \tau_2^1 \circ \tau_3^2
\circ \dots \circ \tau_{i-1}^{i-2} 
\circ \tau_i^{i-1} \circ T^i  \\
\label{another-cup}
\cup_{i,n} & = & T_i \circ \tau_{i-1}^i \circ \tau_{i-2}^{i-1}
\dots \circ \tau_2^3 \circ \tau_1^2 \circ \zeta_{n+2}^{-1} 
\end{eqnarray} 

\begin{conjecture} There are natural isomorphisms
(\ref{first-relation})-(\ref{last-relation}) 
( with  $q$ set to $-1$ in (\ref{last-relation}))
of functors 
where $\cap_{i,n}$ and $\cup_{i,n}$ are defined by
(\ref{another-cap}) and (\ref{another-cup}). 
\end{conjecture}

Equivalences (\ref{first-relation}) and (\ref{second-relation}) 
 are immediate from [ES] and the results of this section.
Relation (\ref{last-relation}) is implied by the fact
 that the composition of translation functors from and to a wall is 
equivalent to two copies of the identity functor. 
 To prove the remaining four equivalences, 
a thorough understanding of the functor $\zeta_n$, 
defined in [ES] in quite a tricky way, will be required. 

To categorify the Temperley-Lieb algebra $TL_{n,q}$ for arbitrary $q$, 
rather than just $q=-1,$ one needs to work with the mixed version of 
parabolic categories and projective functors. The conjecture of 
Irving [Ir] that projective functors admit a  mixed structure 
can probably be approached via a recent work [BGi] of  Beilinson and Ginzburg 
where wall-crossing functors are realized geometrically. 

We next state the parabolic analogue of  
Conjecture~\ref{tangle-invariants-sing}.
 It is convenient to suppress parameter $k$ in the definition 
 of $\cap_{i,n}$ and $\cup_{i,n}$ by summing over $k$ and 
 passing to categories $\O^n= \oplus_k \O^{k,n-k}.$ 
  We switch to derived categories and extend functors $\cap_{i,n}, 
 \cup_{i,n}$ and $\CU_i$ to derived functors 
 \begin{eqnarray*} 
  \cap_{i,n} & : & D^b(\O^n) \lra D^b(O^{n-2}) \\
  \cup_{i,n} & : & D^b(\O^n) \lra D^b(\O^{n+2}) \\
  \CU_i      & : & D^b(\O^n) \lra D^b(\O^n) 
  \end{eqnarray*} 

Recall elementary tangles (a)-(c) described at the end of 
Section~\ref{TL-and-Z}. 

To the 
tangle (a) associate functor $\cap_{i,n}$ given by the formula 
(\ref{another-cap}). 

To the tangle (b) associate functor $\cup_{i,n}$ given by the formula 
(\ref{another-cup}). 

To the tangle (c) associate functor $R_{i,n}$ from $D^b(\O^n)$ to $D^b(\O^n)$ 
which is the cone of the adjointness morphism of functors 
$$ \CU_i  \lra Id. $$
Given a presentation $\alpha$ of a tangle $t$ as a composition 
of elementary tangles of types (a)-(c), to $\alpha$ we associate the 
functor $g(\alpha)$ which is the corresponding composition of 
functors $\cap_{i,n}, \cup_{i,n}, R_{i,n}.$ 

\begin{conjecture}
\label{tangle-invariants-parab} 
 Given two presentations $\alpha, \beta$ of 
a tangle $t$ as products of elementary tangles (a)-(c), 
functors $g(\alpha)$ and $g(\beta)$ are  
isomorphic, up to shifts in the derived category. 
\end{conjecture}

\subsection{$\mf{sl}_2$ and Zuckerman functors} 
\label{sl2-and-Z} 

The Grothendieck group of the category  $\O^n=\oplus_{k=0}^n \O^{k,n-k}$ 
has rank $2^n$. In the previous section we showed that
the projective functors, restricted to $\O^n$, ``categorify'' 
 the Temperley-Lieb 
algebra action on $\vn.$
The Lie algebra $\mf{sl}_2$ action on $\vn$ commutes with the Temperley-Lieb 
algebra action, while Zuckerman functors commute with projective functors. 
It is an obvious guess now that Zuckerman functors between different 
blocks of $\O^n$ provide a ``categorification'' of this $\mf{sl}_2$
action. This fact is well-known and dates back to [BLM] and [GrL], where 
it is presented in a different language
and in the more general case of $\mf{sl}_k$ rather than $\mf{sl}_2.$ 
 Beilinson, Lusztig and 
MacPherson in [BLM] count points over finite fields 
in certain correspondences between flag varieties. These correspondences 
define functors between derived categories of sheaves on these 
flag varieties, smooth along Schubert stratifications. 
Counting points is equivalent to computing the 
action of the corresponding functors on Grothendieck groups. 

These derived categories are equivalent to derived categories of
parabolic subcategories of a  regular block of the highest weight 
category for $\mf{gl}_n.$ Pullback and pushforward functors are 
then isomorphic (up to shifts) to the embedding functor from smaller
to bigger parabolic subcategories and its adjoint 
functors which are Zuckerman functors (see [BGS], Remark (2) on page 504). 
Combining this observation with  the computation 
of [GrL] for the special case of the Grassmannian rather than an arbitrary 
partial flag variety, one can check that Zuckerman functors
between various pieces of $\O^n$ categorify the action of
$E^{(a)}$ and $F^{(a)}$ on $\vn.$ 
This fact is stated more accurately below.

Fix $n\in \N.$ Recall that we denoted by $\mf{p}_k$ the parabolic subalgebra 
of $\mf{gl}_n$ consisting of block upper-triangular matrices with 
the reductive part $\mf{gl}_k\oplus \mf{gl}_{n-k}.$ 
Denote by $\mf{p}_{k,l}$ the parabolic subalgebra of $\mf{gl}_n$ which 
is the intersection of $\mf{p}_k$ and $\mf{p}_{k+l}.$ The maximal 
reductive subalgebra of $\mf{p}_{k,l}$ is isomorphic to 
$\mf{gl}_k\oplus \mf{gl}_{l}\oplus \mf{gl}_{n-k-l}.$ 

Denote by $\O^{k,l,n-k-l}$ the complete subcategory of $\O_{\mu}$ consisting
of $U(\mf{p}_{k,l})$-locally finite modules. We have embeddings of 
categories 
\begin{eqnarray}
\mc{I}_{k,l}: & \O^{k,n-k} & \lra \O^{k,l,n-k-l} \\
\mc{J}_{k,l}: & \O^{k+l,n-k-l} & \lra \O^{k,l,n-k-l} 
\end{eqnarray} 
Denote by $\mc{K}_{k,l}$ and $\mc{L}_{k,l}$ derived functors 
of the right adjoint functors 
of $\mc{I}_{k,l}$ and $\mc{J}_{k,l}$: 

\begin{eqnarray}
\mc{K}_{k,l}: & D^b(\O^{k,l,n-k-l}) & \lra D^b(\O^{k,n-k}) \\
\mc{L}_{k,l}: & D^b(\O^{k,l,n-k-l}) & \lra D^b(\O^{k+l,n-k-l}) 
\end{eqnarray} 
These are derived functors of Zuckerman functors. Compositions 
$\mc{L}_{k,l}\circ \mc{I}_{k,l}$ and $\mc{K}_{k,l}\circ 
\mc{J}_{k,l}$ are exact functors between derived categories 
$D^b(\O^{k,n-k})$ and $D^b(\O^{k+l,n-k-l}),$ and  in the Grothendieck 
group of $\O^n=\oplus_{k=0}^n \O^{k,n-k}$ 
descend to the action of $E^{(l)}$ and $F^{(l)}$ on 
the module $\vn.$ 

\hspace{0.5in}

\end{document}